\documentclass[11pt,reqno]{amsart}
\usepackage[a4paper,margin=30mm,headheight=14pt]{geometry}
\usepackage{amsmath,amssymb,amsthm,mathtools,mathrsfs,booktabs,array}
\usepackage{xcolor,enumitem,etoolbox,microtype,needspace}
\usepackage[T1]{fontenc}
\usepackage{lmodern}
\usepackage{hyperref}
\usepackage{aliascnt}
\usepackage[nameinlink,capitalize,noabbrev]{cleveref}
\hypersetup{
  colorlinks=true,
  linkcolor=blue!45!black,
  citecolor=blue!45!black,
  urlcolor=blue!45!black
}
\hypersetup{
  pdftitle={Global rigidity and the Hessian energy identity for the one-phase Bernoulli problem},
  pdfauthor={Zhouyu Long and Wenming Zou},
  pdfsubject={Classical stable solutions and bounded-index curvature concentration},
  pdfkeywords={Bernoulli problem, free boundary, stability, Morse index, Hessian energy}
}
\numberwithin{equation}{section}

\newtheorem{theorem}{Theorem}[section]
\theoremstyle{plain}
\newaliascnt{lemma}{theorem}
\newtheorem{lemma}[lemma]{Lemma}
\aliascntresetthe{lemma}
\theoremstyle{plain}
\newaliascnt{proposition}{theorem}
\newtheorem{proposition}[proposition]{Proposition}
\aliascntresetthe{proposition}
\theoremstyle{plain}
\newaliascnt{corollary}{theorem}
\newtheorem{corollary}[corollary]{Corollary}
\aliascntresetthe{corollary}
\theoremstyle{definition}
\newaliascnt{definition}{theorem}

\aliascntresetthe{definition}
\theoremstyle{remark}
\newaliascnt{remark}{theorem}
\newtheorem{remark}[remark]{Remark}
\aliascntresetthe{remark}
\theoremstyle{remark}
\newaliascnt{example}{theorem}

\aliascntresetthe{example}
\theoremstyle{plain}
\crefname{theorem}{Theorem}{Theorems}
\Crefname{theorem}{Theorem}{Theorems}
\crefname{lemma}{Lemma}{Lemmas}
\Crefname{lemma}{Lemma}{Lemmas}
\crefname{proposition}{Proposition}{Propositions}
\Crefname{proposition}{Proposition}{Propositions}
\crefname{corollary}{Corollary}{Corollaries}
\Crefname{corollary}{Corollary}{Corollaries}
\crefname{definition}{Definition}{Definitions}
\Crefname{definition}{Definition}{Definitions}
\crefname{remark}{Remark}{Remarks}
\Crefname{remark}{Remark}{Remarks}
\crefname{example}{Example}{Examples}
\Crefname{example}{Example}{Examples}
\newcommand{\R}{\mathbb R}
\newcommand{\one}{\mathbf 1}
\newcommand{\diver}{\operatorname{div}}
\newcommand{\dd}{\,\mathrm d}
\newcommand{\FB}{\operatorname{FB}}
\newcommand{\Exc}{\mathcal E}

\title[Global rigidity and the Hessian energy identity]{\textrm{Global rigidity and the Hessian energy identity\\ for the one-phase Bernoulli problem}}
\author[Zhouyu Long]{\bfseries Zhouyu Long\(^{*}\)}
\author[Junbao Zhou]{\bfseries  Junbao Zhou\(^{\dagger}\)}
\thanks{\(^{*}\)Department of Mathematical Sciences, Tsinghua University,
Beijing 100084, China.\newline
E-mail: \href{mailto:longzy25@mails.tsinghua.edu.cn}{\nolinkurl{longzy25@mails.tsinghua.edu.cn}}.}
\thanks{\(^{\dagger}\)Department of Mathematical Sciences, Tsinghua University,
Beijing 100084, China.\newline
E-mail: \href{mailto:zjb23@mails.tsinghua.edu.cn}{\nolinkurl{zjb23@mails.tsinghua.edu.cn}}.}

\makeatletter
\patchcmd{\@setauthors}{\centering\footnotesize}{\centering\large}{}%
  {\ClassError{amsart}{Unable to set author font}{Check the amsart author layout.}}
\patchcmd{\@setauthors}{\MakeUppercase{\authors}}{\authors}{}%
  {\ClassError{amsart}{Unable to preserve author case}{Check the amsart author layout.}}
\def\author@andify{\nxandlist{\unskip,\quad}{\unskip,\quad}{\unskip,\quad}}
\let\addresses\@empty
\makeatother
\date{}

\begin{document}

\begin{abstract}
We prove that every entire classical stable solution of the one-phase
Bernoulli problem in dimension {\bf five or six}  is one-dimensional.    In dimensions
three through six, we establish an exact interior identity for the
scale-invariant Hessian energy of sequences with bounded Morse index. 
Hence, the global stability conjecture is solved. 

\end{abstract}

\maketitle



\section{\bf Introduction}\label{sec:intro}

Let \(D\subset\R^n\) be an open domain, and consider a continuous nonnegative
function \(u:D\to[0,\infty)\) with positivity set (positive phase)
\(\Omega=\{x\in D:u(x)>0\}\) and free boundary
\(\FB(u)=\partial\Omega\cap D\). We say that \(u\) is a \emph{classical
solution} to the one-phase Bernoulli free boundary problem in \(D\) if
\(\FB(u)\) is a smooth hypersurface around which \(\Omega\) is locally a
smooth one-sided domain, and \(u\) satisfies the overdetermined boundary value
problem
\begin{equation}\label{eq:bernoulli}
  \Delta u=0\quad\text{in }\Omega,\qquad
  u=0,\quad |\nabla u|=1\quad\text{on }\FB(u),
\end{equation}
where \(\nabla u\) along \(\FB(u)\) denotes the non-tangential trace from
\(\Omega\). Introduced in its modern variational formulation by Alt and
Caffarelli \cite{AC}, \eqref{eq:bernoulli} describes the Euler--Lagrange equilibrium of
the functional
\[
  J(u;D)=\int_D \bigl(|\nabla u|^2+\one_\Omega\bigr)\dd x
\]
under compactly supported inner domain variations. Denoting by
\(\nu=-\nabla u\) the outward unit normal to \(\Omega\) along \(\FB(u)\) and
by \(H=-\diver_{\FB(u)}\nu\) the corresponding scalar mean curvature (signed
so that \(H=n-1>0\) on the unit sphere bounding the exterior-ball domain
\(\R^n\setminus\overline{B_1}\)), we use the scalar second-variation form of
\cite[Lemma~1]{CJK} in the convention of \cite[(1.2)]{FFSS}
\[
  \mathcal Q_u[\phi;D]
  =\int_\Omega |\nabla\phi|^2\dd x
   -\int_{\FB(u)} H\phi^2\dd\mathcal H^{n-1}.
\]
Here admissible test functions \(\phi\) are restrictions to
\(\overline\Omega\cap D\) of smooth functions in \(C_c^\infty(D)\), which need
not vanish on \(\FB(u)\). The \emph{Morse index} \(\operatorname{ind}(u;D)\) is
defined as the supremum of the dimensions of subspaces of admissible test
functions on which \(\mathcal Q_u[\cdot;D]\) is negative definite; when
\(D=\R^n\), we simply write \(\operatorname{ind}(u)\). In particular, \(u\) is
\emph{stable} in \(D\) if \(\operatorname{ind}(u;D)=0\), that is,
\(\mathcal Q_u[\phi;D]\ge0\) for all admissible \(\phi\).

The one-phase problem has antecedents in free-streamline models in fluid
mechanics. The periodic hollow vortices constructed by Baker, Saffman and
Sheffield \cite{BSS} provide a classical example of this connection.
Alt and Caffarelli \cite{AC} introduced the modern variational framework and
established existence and regularity for minimizers. Caffarelli's
Harnack-inequality approach subsequently established regularity of Lipschitz
free boundaries and a flatness-to-Lipschitz theorem
\cite{CaffarelliI,CaffarelliII}. We refer to \cite{CS} for the geometric
theory of these problems. Weiss \cite{WeissWeak,WeissMin} developed the
monotonicity formula and used one-homogeneous minimizing blow-ups to relate
the dimension of the singular set to the first dimension admitting a
nonflat minimizing cone.

For one-homogeneous stable solutions with free boundary smooth away from
the vertex, the rigidity results of Caffarelli--Jerison--Kenig \cite{CJK}
and Jerison--Savin \cite{JS} give flatness through dimensions three and four,
respectively. The classification of these smooth stable cones in dimensions
five and six follows from the analysis underlying the minimizing-cone
theorem in \cite[Sections~3--7]{OAI}, as recorded in
\cite[Theorem~3 and Remark~4]{survey}. In dimension seven,
De Silva and Jerison \cite{DSJ} constructed a nonflat minimizing cone.
The foliation theorem of Edelen, Spolaor and Velichkov
\cite[Theorem~1.5]{ESV} gives smooth entire minimizing leaves with that
cone as a blow-down.

A central problem in free boundary regularity---parallel to the Bernstein
problem for minimal hypersurfaces and De Giorgi's conjecture for the
Allen--Cahn equation---is the \emph{global stability conjecture} for
\eqref{eq:bernoulli} (see \cite[Conjecture~6]{survey}):
\begin{quote}
\itshape
Up to which ambient dimension must every entire classical stable solution
\(u:\R^n\to[0,\infty)\) be one-dimensional?
\end{quote}

In dimension \(n=2\), rigidity follows from a logarithmic cutoff
argument \cite[discussion following Conjecture~6]{survey}.
Chan, Fern\'andez-Real, Figalli and Serra \cite{CFFS} resolved the conjecture
in dimension \(n=3\), and Fern\'andez-Real and Serra \cite{FRS} established
the case \(n=4\) using a multiscale signed decomposition adapted to
non-minimizing two-sheeted configurations. The resulting gap between the
smooth stable cone classification and the entire stable classification
is in dimensions five and six.

\subsection{Main results}\label{subsec:main-results}

Our first objective is to close this dimensional gap completely. To separate
our analytic PDE reduction cleanly from the external cone-classification input,
recall following \cite[Definition~4.1]{FRS} that the critical cone dimension
\(n_{\mathrm{cone}}^*\) is the smallest integer \(n\) admitting a nonzero,
non-half-space one-homogeneous global stable solution that is classical away
from the origin with \(0\) on its free boundary; similarly, let
\(n_{\mathrm{entire}}^*\) denote the smallest dimension admitting a
non-one-dimensional entire classical stable solution. Our first main theorem
shows that cone rigidity transfers unconditionally to entire rigidity in all
dimensions below \(n_{\mathrm{cone}}^*\).

\begin{theorem}\label{thm:transfer}
If \(4\le n<n_{\mathrm{cone}}^*\), then every entire classical stable solution
\(u\) to \eqref{eq:bernoulli} in \(\R^n\) satisfies \(D^2u\equiv0\) in its
positive phase \(\Omega\) and is one-dimensional.
\end{theorem}

The stable-cone classification through dimension six
\cite[Theorem~3 and Remark~4]{survey} and the seven-dimensional
De Silva--Jerison cone \cite{DSJ} give \(n_{\mathrm{cone}}^*=7\). Combining Theorem~\ref{thm:transfer} with the two- and
three-dimensional results \cite{CFFS,survey} yields the full classification of
entire classical stable solutions through dimension six.

\begin{corollary}\label{cor:classification}
Let \(2\le n\le6\). Every entire classical stable solution \(u\) to
\eqref{eq:bernoulli} in \(\R^n\) is either a nonnegative constant or, for some
unit vector \(e\in\mathbb S^{n-1}\), one of the one-dimensional half-space or
two-plane profiles
\[
  (x\cdot e-a)_+\quad (a\in\R),
  \qquad\text{or}\qquad
  (x\cdot e-a)_++(b-x\cdot e)_+\quad (b<a).
\]
In particular, \(D^2u\equiv0\) in \(\Omega\), and \(n_{\mathrm{entire}}^*=7\)
is the sharp first dimension admitting a nonflat entire classical stable
solution (see Corollary~\ref{cor:critical} below).
\end{corollary}

By a point-picking blow-up argument, entire rigidity immediately implies
universal, scale-invariant interior bounds on both the positive-phase Hessian
and the second fundamental form \(\mathrm{II}\) of the free boundary.
The connection between entire rigidity and local curvature estimates is
already present in the conditional estimate of Kamburov and Wang \cite{KW}.

\begin{corollary}\label{cor:curv}
Let \(2\le n\le6\). There exists a dimensional constant \(C_n>0\) such that if
\(u\) is a classical stable solution to \eqref{eq:bernoulli} in
\(B_R(z)\subset\R^n\) with \(z\in\FB(u)\), then
\[
  \|D^2u\|_{L^\infty(B_{R/2}(z)\cap\{u>0\})}
  +\|\mathrm{II}\|_{L^\infty(B_{R/2}(z)\cap\FB(u))}
  \le C_n R^{-1},
\]
where \(D^2u\) extends continuously up to \(\FB(u)\) from the positive phase.
\end{corollary}

For planar solutions with connected positivity set, Basulto and Kamburov
\cite[Theorem~1.1]{BK} classified the finite-index profiles. In higher
dimensions, Liu, Wang and Wei \cite{LWW} constructed axially symmetric
catenoid-type solutions. The cone classification and
Corollary~\ref{cor:classification} give, respectively,
\(n_{\mathrm{cone}}^*>6\) and \(n_{\mathrm{entire}}^*>6\).
Together these verify the threshold required by the recent work of Fern\'andez-Real, Florit-Simon and Serra
\cite[Theorem~1.2 and (1.3)--(1.5)]{FFSS}. Consequently, every non-flat entire
classical solution of finite Morse index in dimensions \(3\le n\le6\) is
axially symmetric and belongs, up to rigid motion and Bernoulli rescaling
\(u_{z,\lambda}(x)=\lambda^{-1}u(z+\lambda x)\), either to the radial
exterior-ball family or to the two-ended \emph{catenoid-type} family (recorded
precisely in Corollary~\ref{cor:finite}).

Once all entire finite-Morse-index blow-up models are classified, it is
natural to ask how curvature concentrates along sequences of classical
solutions with uniformly bounded Morse index. A geometric precedent is the
total-curvature quantization theorem of Buzano and Sharp \cite{BS} for closed
minimal hypersurfaces of bounded area and index. In the Bernoulli problem,
the scale-invariant curvature measure we study is the \(L^n\) Hessian measure
\[
  \mu_u(E)=\int_{E\cap\{u>0\}}|D^2u|^n\dd x,
\]
where \(|D^2u|\) denotes the Hilbert--Schmidt norm of the classical Hessian in
\(\{u>0\}\) (extended by zero outside \(\{u>0\}\)), which satisfies
\(\mu_{u_{z,\lambda}}(E)=\mu_u(z+\lambda E)\). In the Bernoulli setting,
classical free boundaries can collapse onto two-sided
\emph{unit wedge} limits \(|(x-p)\cdot e|\) whose free boundaries develop
self-contact without carrying positive-phase Hessian mass, while shrinking
exterior-ball solutions converge uniformly to zero even as their positivity
sets fill the ambient ball in \(L^1\). Kriventsov and Weiss
\cite[Theorems~1.2 and~9.3]{KriventsovWeiss} developed a compact variational
formulation that tracks the phase indicator independently of the function
in such limits. Here the second main theorem establishes an exact interior
bubble-tree identity for the Hessian measure. The proof controls its mass
across intermediate neck scales.

\begin{theorem}\label{thm:local-bubble-tree}
Let \(3\le n\le6\), and let \(\{u_j\}\) be a sequence of classical solutions
to \eqref{eq:bernoulli} in \(B_2\subset\R^n\) satisfying
\(\operatorname{ind}(u_j;B_2)\le I\) and \(\FB(u_j)\cap B_1\ne\varnothing\).
After passing to a subsequence, \(u_j\) converges locally uniformly in \(B_2\)
to a locally Lipschitz function \(U\) that is harmonic on \(\{U>0\}\) (whose
free boundary may exhibit two-sided sheet contact), and there exist an integer
\(0\le N\le I\), nonflat entire classical solutions \(V^1,\ldots,V^N\) of
finite Morse index, centers \(z_j^i\in\FB(u_j)\) converging to points
\(p^i\in B_2\), and blow-up scales \(\lambda_j^i\to0\) with the following
properties:
\begin{enumerate}
\item[\textup{(i)}] For each \(1\le i\le N\), the rescalings
\(u_j(z_j^i+\lambda_j^i x)/\lambda_j^i\) converge locally uniformly on
\(\R^n\) to \(V^i\), smoothly on compact subsets of the positive phase in
local free-boundary graph charts, and with \(\one_{\{u_j(z_j^i+\lambda_j^i\cdot)>0\}}\to\one_{\{V^i>0\}}\)
in \(L^1_{\mathrm{loc}}(\R^n)\). Each bubble \(V^i\) is, up to rigid motion
and Bernoulli scaling, either the radial exterior-ball solution or an axially
symmetric catenoid-type solution, and has finite total Hessian energy
\(\mu_{V^i}(\R^n)<\infty\).
\item[\textup{(ii)}] The bubble scales and centers are pairwise separated in
the bubble-tree sense, and the total Morse index is bounded by \(I\):
\[
  \sum_{i=1}^N\operatorname{ind}(V^i)\le I,\qquad
  \frac{|z_j^i-z_j^k|}{\max(\lambda_j^i,\lambda_j^k)}\longrightarrow\infty
  \quad (1\le i\ne k\le N).
\]
\item[\textup{(iii)}] As Radon measures locally in \(B_2\), the
scale-invariant Hessian measures satisfy the exact energy identity
\begin{equation}\label{eq:local-multi-energy}
  \mu_{u_j}\rightharpoonup
  |D^2U|^n\one_{\{U>0\}}\dd x
  +\sum_{i=1}^N \mu_{V^i}(\R^n)\delta_{p^i},
\end{equation}
where the concentration points \(p^1,\ldots,p^N\) need not be distinct.
\end{enumerate}
\end{theorem}

As immediate consequences of Theorem~\ref{thm:local-bubble-tree}, we obtain
uniform interior \(L^n\)-Hessian bounds depending only on \(n\), \(I\), and
the subdomain (see Corollary~\ref{cor:local-uniform-energy} below). When curvature blows
up at \(I\) distinct points, there is one index-one bubble at each point
(see Corollary~\ref{cor:local-index-saturated} below). Proposition~\ref{prop:exterior-ball-energy}
computes the exterior-ball bubble energy
$$E_n^{\mathrm{ball}}=|\mathbb S^{n-1}|[n(n-1)]^{n/2-1}.$$
For a two-ended bubble, the weight in \eqref{eq:local-multi-energy} retains
its dependence on the profile. Proposition~\ref{prop:seven-energy} gives
index-zero rescalings in dimension seven for which the interior Hessian
energy grows at least logarithmically. Thus the dimensional
threshold \(n\le6\) is sharp for both parts of the paper.

\subsection{Strategy of the proofs and organization of the paper}\label{subsec:strategy}

Let us outline the main analytical mechanisms behind Theorems~\ref{thm:transfer}
and~\ref{thm:local-bubble-tree}.

\smallskip
\noindent\emph{1. Localized signed tree and fixed-center excess contraction.}
Suppose toward a contradiction that a nonflat entire classical stable solution
exists in dimension \(4\le n<n_{\mathrm{cone}}^*\). After normalizing
\(\|\nabla u\|_{L^\infty}\le1\) and \(\|D^2u\|_{L^\infty}=|D^2u(0)|=1\), the
\(\varepsilon\)-regularity theory and cone classification below
\(n_{\mathrm{cone}}^*\) produce a nonempty set of intrinsic neck centers
\(\mathcal Z\) with regularity scales \(r_*(\zeta)\ge c_0>0\) carrying a
definite quantum \(\eta_0^n\) of \(L^n\)-Hessian mass. This setup comes
from \cite{CFFS} and is developed in \cite[Section~5]{FRS}.
A major obstacle when passing from \(n=4\) to \(n=5,6\) is that the
\((n-2)\)-dimensional neck content \(Q(r)=\sum\ell_N^{n-2}\) has supercritical
scaling \(r^{n-2}\) at infinity. Rather than performing a global blow-down of
trees across \(\R^n\), we fix a single neck center \(z\in\mathcal Z\) and build
a finite Calderón--Zygmund-type signed stopping-time tree inside one large
root ball \(B_R(z)\), decomposing \(\Omega\cap B_{R/2}(z)\) outside a union of
neck terminal balls into two disjoint Caccioppoli sheets \(U_\pm\) closely
modeled on a symmetric unit wedge \(|(x-z)\cdot e|\)
(Lemma~\ref{lem:fixed-root-tree}). Because \(z\in\mathcal Z\) is fixed, the
neck terminal \(N_*\) containing \(z\) has a strictly positive radius
\(\ell_*>0\) that persists at every scale \(r>0\) around \(z\).

To bound the neck content \(Q(s)\) as we zoom in toward \(z\), we project the
gradient \(\nabla u\) onto an \((n-2)\)-plane \(M\) nearly orthogonal to both
approximate sheet normals. The components of \(F=P_M\nabla u\) are translation
Jacobi fields, so they obey a Sternberg--Zumbrun-type stability inequality
\cite{SZ}, in the vector-valued Bernoulli formulation of
\cite[Lemma~6.3]{FRS}. This controls the directional curvature density
\(\mathcal A_M^2=|F|^2|D(F/|F|)|^2\).
In Lemma~\ref{lem:gap}, we prove a uniform coercivity gap on each doubled
intrinsic ball \(B_{2r_*(\zeta)}(\zeta)\): if \(\mathcal A_M^2\) vanished in
the blow-up limit, harmonic unique continuation would produce \(n-3\)
independent translation-invariant directions inside \(M\), splitting off an
\(\R^{n-3}\) factor and yielding a classical stable quotient in \(\R^3\),
which must have vanishing Hessian by the three-dimensional theorem of
\cite{CFFS}. Summing this local gap against the Sternberg--Zumbrun cutoff
estimate bounds \(Q(\theta r)\) quadratically in terms of the projected
Dirichlet energy. In Lemma~\ref{lem:contract}, we couple the two-sided
\(L^2\) affine excess on \(U_\pm\cap B_r\) with \(\kappa r^{2-n}Q(r)\): a
divergence-free Rellich stress-energy identity controls free-boundary normal
deviation and rules out \(L^2\) concentration in the vanishing nodal slab via
the Euclidean \(BV\)--Sobolev inequality, allowing harmonic linearization with
homogeneous Neumann data across the nodal hyperplane. Iterating the resulting
excess contraction \(\Exc(\theta^j r_0,L_j)\le 2^{-j}\Exc(r_0,L_0)\) on a
single sufficiently large tree forces \(Q(\theta^j r_0)\to0\), contradicting
\(Q(\theta^j r_0)\ge\ell_*^{n-2}>0\).

\smallskip
\noindent\emph{2. Weiss density pinning, annular splitting, and Voronoi bubble trees.}
For a sequence \(\{u_j\}\) with \(\operatorname{ind}(u_j;B_2)\le I\),
repeatedly applying the weighted point-picking argument of
Lemma~\ref{lem:weighted-selection} extracts at most \(I\) scale-separated
bubbles \(V^1,\ldots,V^N\) and leaves a pointwise curvature bound inversely
proportional to the distance \(d_j(x)\) to the bubble centers. The central
difficulty in proving the exact measure identity \eqref{eq:local-multi-energy}
is to rule out diffuse Hessian mass across intermediate annuli
\(\lambda_j^i\ll s_j\ll 1\) and across interaction scales
\(|z_j^i-z_j^k|\) when several bubbles converge to the same point \(p\).

Our key observation (i.e., Lemma~\ref{lem:weiss-cap} and
Proposition~\ref{prop:intermediate-host}) is that for both possible bubble
types---radial exterior balls and axially symmetric catenoid-type
solutions---the asymptotic Weiss energy at infinity equals the full unit-ball
volume \(\omega_n=|B_1|\). On the other hand, integrating the Weiss
monotonicity differential inequality from any interior scale up to a fixed
macroscopic ball caps \(W(u_j,r)\) from above by \(\omega_n+O(1/\log(R/r))\).
Sandwiching the Weiss functional between the bubble scale and the macroscopic
scale pins \(W(w_j,t)\to\omega_n\) at every intermediate scale \(s_j\). This
simultaneously forces one-homogeneity of any intermediate limit \(w\) and pins
the weak-\(*\) limit of the positivity indicators to \(\chi\equiv1\) almost
everywhere. Outside at most \(I\) unstable concentration points, our local
curvature estimate (Corollary~\ref{cor:curv}) and local multi-sheet analysis
(Lemmas~\ref{lem:local-sheets}--\ref{lem:sheet-contact}) then rule out
one-sided sheets (which would leave an open zero phase contradicting
\(\chi\equiv1\)) and restrict \(w\) to be either \(0\) or a centered unit
wedge \(|x\cdot e|\). Across each scale gap, a topological connectedness
dichotomy (Lemmas~\ref{lem:glue-neck-branches}--\ref{lem:abstract-neck}) shows
that \(u_j\) is either globally harmonic on the gap annulus---where multipole
harmonic decay applies---or splits globally into two disjoint classical
branches \(v_\pm\) to which the annular improvement-of-flatness theorem of
\cite{FFSS} applies with polynomial decay rate \(1/2\). Finally, partitioning
space into Voronoi nearest-center cells \(E_j^i\) around colliding centers
isolates every transition scale away from neighboring singularities, completing
the no-neck energy accounting.

\smallskip
\noindent\textbf{Organization of the paper.}
Section~\ref{sec:cone-transfer} establishes the localized signed tree
decomposition, the projected Sternberg--Zumbrun neck estimate, and the
fixed-center excess contraction, proving Theorem~\ref{thm:transfer},
Corollary~\ref{cor:classification}, and the sharp seven-dimensional threshold
(Corollary~\ref{cor:critical}). Section~\ref{sec:curvature} develops local
gradient bounds, multi-sheet graphical compactness, and the point-picking
lemma, deducing Corollary~\ref{cor:curv} and the finite-index classification
(Corollary~\ref{cor:finite}). Section~\ref{sec:index-concentration} proves the
logarithmic Weiss bound, intermediate-scale density pinning, annular energy
vanishing, and Theorem~\ref{thm:local-bubble-tree}. Finally,
Section~\ref{sec:examples} computes the exterior-ball energy and
constructs the seven-dimensional example with logarithmic lower growth.

\section{\bf Global rigidity below the critical cone dimension}\label{sec:cone-transfer}

Throughout Sections~\ref{sec:tree}--\ref{sec:contraction}, we argue by
contradiction to prove Theorem~\ref{thm:transfer}. Suppose that in some
dimension \(4\le n<n_{\mathrm{cone}}^*\) there exists a nonflat entire
classical stable solution to \eqref{eq:bernoulli}. By the global gradient bound
and normalization procedure of \cite[Lemmas~3.5 and~5.2]{FRS}, we may pass to
an entire classical stable solution \(u\) satisfying \(|\nabla u|\le1\) and
\(|D^2u|\le1\) everywhere in \(\Omega\), together with \(0\in\FB(u)\) and
\(|D^2u(0)|=1\) (understood as the norm of the positive-side trace).

\vskip0.13in

Following \cite[Definition~5.4]{FRS}, let \(\mathcal Z\subset\FB(u)\)
denote the collection of intrinsic neck centers of \(u\), with associated
regularity scales \(r_*(\zeta)\) characterized by the threshold Hessian mass
\[
  \int_{\Omega\cap B_{r_*(\zeta)}(\zeta)}|D^2u|^n\dd x=\eta_0^n>0,
  \qquad \inf_{\zeta\in\mathcal Z}r_*(\zeta)>0.
\]
Because \(|D^2u(0)|=1\) and \(\|D^2u\|_{L^\infty(\Omega)}\le1\), small-Hessian
\(\varepsilon\)-regularity ensures that \(\mathcal Z\) is nonempty and that
\(r_*(\zeta)\) has a uniform positive lower bound \cite[Section~5.3]{FRS}.
Moreover, all local neck and wedge approximation properties established in
\cite[Section~5]{FRS} hold verbatim under the dimensional hypothesis
\(n<n_{\mathrm{cone}}^*\) (see \cite[Remark~5.1 and Section~5.4]{FRS}).

\subsection{Localized signed tree decomposition at a fixed neck center}\label{sec:tree}

Fix once and for all a neck center \(z\in\mathcal Z\), and abbreviate
\(B_r=B_r(z)\). Rather than letting the tree root vary dynamically, we verify
that for every sufficiently large root radius \(R\), the dyadic stopping-time
construction of \cite[Section~9]{FRS} initialized at \(B_R\) produces two
disjoint Caccioppoli sheets \(U_\pm\subset\Omega\cap B_{R/2}\) whose boundary
defects are controlled by the \((n-2)\)-dimensional neck content, while
retaining a persistent terminal ball around \(z\).

\begin{lemma}\label{lem:fixed-root-tree}
Fix \(z\in\mathcal Z\) and write \(B_r=B_r(z)\). For a sufficiently small fixed
aperture parameter \(\delta>0\) and every sufficiently large root radius \(R\),
there exists a finite Calderón--Zygmund-type signed tree rooted at \(B_R\)
yielding disjoint open sets \(U_+,U_-\subset\Omega\cap B_{R/2}\) of locally
finite perimeter and a finite family of neck terminal balls
\(N=B_{\ell_N}(x_N)\), one of which, denoted \(N_*\), contains \(z\). With
constants independent of \(R\), the decomposition satisfies:
\begin{enumerate}
\item[\textup{(i)}] Each neck terminal \(N\) admits an associated center
\(\zeta_N\in\mathcal Z\cap 4N\) with comparable scale
\(r_*(\zeta_N)\asymp_\delta\ell_N\). The map \(N\mapsto\zeta_N\) has bounded
multiplicity, the doubled intrinsic balls \(B_{2r_*(\zeta)}(\zeta)\) of
distinct centers in \(\mathcal Z\) have bounded overlap, and
\(\max_N\ell_N=o(R)\) as \(R\to\infty\).

\item[\textup{(ii)}] Whenever \(0<r\le R/32\) and the scale ratio \(r_*(z)/r\)
is sufficiently small, there exists a unit vector \(e\in\mathbb S^{n-1}\) such
that, in terms of the signed normal coordinate \(t=(x-z)\cdot e\),
\begin{equation}\label{eq:slab}
  B_r\cap\{\pm t>\delta_r r\}\subset U_\pm\cap B_r
  \subset B_r\cap\{\pm t>-\delta_r r\},\qquad
  \|u-|t|\|_{L^\infty(B_r)}=o(r),
\end{equation}
where both \(\delta_r\) and \(r^{-1}\|u-|t|\|_{L^\infty(B_r)}\) tend to zero
as \(r_*(z)/r\to0\), uniformly in \(R\).

\item[\textup{(iii)}] Let \(Y=\Omega\setminus(U_+\cup U_-)\) denote the
residual positive phase and \(G_\pm=\partial^*U_\pm\cap\Omega\) the interior
reduced boundaries of the signed sheets. Inside \(B_{R/8}\), both \(Y\) and
\(G_\pm\) are covered by the union of closed neck terminals \(\overline N\).
Defining for \(0<r\le R/32\) the \((n-2)\)-dimensional neck content and
maximal neck radius
\[
  Q(r)=\sum_{\overline N\cap B_r\ne\varnothing}\ell_N^{n-2},
  \qquad \ell(r)=\max_{\overline N\cap B_r\ne\varnothing}\ell_N,
\]
the index set is nonempty (as \(N_*\cap B_r\ne\varnothing\)), and one has the
perimeter and volume bounds
\begin{equation}\label{eq:area-volume}
  S(r):=\sum_\pm\mathcal H^{n-1}(G_\pm\cap B_r)\le C\ell(r)Q(r),
  \qquad |Y\cap B_r|\le C\ell(r)^2Q(r).
\end{equation}
\end{enumerate}
\end{lemma}

\begin{proof}
We divide the verification into four steps.

\smallskip
\noindent\emph{Step 1: Stopping-time tree initialization and propagation of intrinsic scales.}
We run the dyadic free-boundary net, local-child generation, and stopping
rules of \cite[Definitions~9.1 and~9.3]{FRS} starting from the root ball
\(B_R=B_R(z)\). For any active ball \(B=B_{r_B}(x_B)\), let \(h_B\) denote its
normalized \(L^2\) wedge excess (measuring \(L^2\) proximity to symmetric
unit wedges on \(16B\), as in \cite[(9.2)]{FRS}). Recall that \(B\) branches
precisely when \(4B\cap\mathcal Z\ne\varnothing\) and \(h_B<\delta\); it stops
as a \emph{regular terminal} if \(4B\cap\mathcal Z=\varnothing\), and as a
\emph{neck terminal} if \(4B\cap\mathcal Z\ne\varnothing\) while
\(h_B\ge\delta\). Because \(z\in\mathcal Z\) is fixed, \cite[Lemma~5.6]{FRS}
implies \(h_{B_R}\to0\) as \(R\to\infty\), so the root ball branches for all
large \(R\). Recentering a parent's optimal wedge at each child's
free-boundary center yields the inherited flatness bound
\(h_B\le C_{\mathrm{inh}}\delta\) for every active ball in the tree
\cite[Lemma~9.4]{FRS}.
\vskip0.12in

Fix \(\rho\in(0,1/16]\). By the conditional step of \cite[Lemma~9.9]{FRS},
whenever an active ball \(B\) and a center \(\zeta\in\mathcal Z\cap 8B\)
satisfy \(r_*(\zeta)\le\rho r_B\), one has the refined upper bound
\begin{equation}\label{eq:tree-small-radii}
  r_*(\zeta)\le C_{\mathrm{nr}}h_Br_B.
\end{equation}
(Indeed, the one-sided barrier of \cite[Corollary~9.2]{FRS} permits recentering
a lower wedge barrier at \(\zeta\) on \(B_{4r_*(\zeta)}(\zeta)\subset 12B\),
where \cite[Lemma~9.8]{FRS} applies.) At the root \(B_R\), no center
\(\zeta\in\mathcal Z\cap B_{8R}\) can have \(r_*(\zeta)\ge\rho R\) once \(R\)
is large: otherwise, scale comparison \cite[Lemma~5.7]{FRS} based at \(\zeta\)
would force \(r_*(z)\ge r_*(\zeta)/C_{8/\rho}\ge \rho R/C_{8/\rho}\),
contradicting the finiteness of \(r_*(z)\). Thus \eqref{eq:tree-small-radii}
holds at the root. Inductively, if \(B'\) is a child of a branching ball \(B\),
then \(8B'\subset(21/4)B\) and \(r_B=2r_{B'}\); since \(h_B<\delta\),
\eqref{eq:tree-small-radii} at \(B\) gives
\(r_*(\zeta)\le 2C_{\mathrm{nr}}\delta r_{B'}\) for all
\(\zeta\in\mathcal Z\cap 8B'\). Choosing \(\delta>0\) small enough that
\(2C_{\mathrm{nr}}\delta\le\rho\) closes the induction and establishes
\eqref{eq:tree-small-radii} across the entire tree. Because every branching
ball meets \(\mathcal Z\) in \(4B\) and \(\inf_{\mathcal Z}r_*>0\), branching
radii are bounded away from zero, and since each dyadic generation has finitely
many net points, the tree terminates in finitely many steps.

\vskip0.13in

\smallskip
\noindent\emph{Step 2: Sheet orientations and terminal scale comparability.}
We propagate normal orientations from parents to children via
\cite[Lemma~9.6]{FRS}. To check same-generation compatibility, let \(D,E\) be
two active balls of equal radius \(s\) with \(8D\cap 8E\ne\varnothing\). Their
parent centers lie at mutual distance strictly less than \(16s+5s=21s\), so the
eightfold dilations of their parents already intersect. The root has a single chosen orientation. At the next generation, compare
the two absolute affine representatives through their \(L^2\) errors with
\(u\) on a ball of radius comparable to \(s\) in \(16D\cap16E\).
Such a ball can be centered at the midpoint of the two centers. Applying
\cite[Lemma~9.5]{FRS} on that ball gives an \(O(\delta)\) gradient difference
and an \(O(\delta s)\) affine difference on \(8D\cap8E\), up to a sign.
The induction hypothesis on the parents and the parent--child gradient
bounds rule out the negative sign. This gives the same uniform bounds at
the new generation as in \cite[Lemma~9.11(ii)]{FRS}.
With this compatibility in place, the absence of holes
\cite[Lemma~9.11(i)]{FRS}, the two-sheet splitting over regular terminals
\cite[Lemma~9.10]{FRS}, and the assembly of signed cores
\cite[Definition~9.12 and Lemma~9.13]{FRS} rely purely on local ball-net
geometry and therefore carry over unchanged, producing disjoint open sets
\(U_\pm\subset\Omega\cap B_{R/2}\) that cover \(\Omega\cap B_{R/2}\) away from
the closed neck terminals.

\vskip0.13in

Now let \(N=B_{\ell_N}(x_N)\) be any neck terminal and pick
\(\zeta_N\in\mathcal Z\cap 4N\). By \eqref{eq:tree-small-radii} and the
inherited bound \(h_N\le C_{\mathrm{inh}}\delta\), we have
\(r_*(\zeta_N)\le C_\delta\ell_N\). Conversely, if \(\ell_N/r_*(\zeta_N)\)
were arbitrarily large relative to \(\delta\), then \cite[Lemma~5.6]{FRS}
would render \(u\) uniformly \(o(\ell_N)\)-close to a centered unit wedge on
\(B_{32\ell_N}(\zeta_N)\); recentering at \(x_N\in\FB(u)\) would yield a
competitor on \(16N\) with \(h_N<\delta\), violating the neck stopping rule.
Hence \(r_*(\zeta_N)\asymp_\delta\ell_N\). Since \(\zeta_N\in 4N\subset B_{5R}\),
applying \eqref{eq:tree-small-radii} at the root gives
\(\ell_N\le C_\delta h_{B_R}R=o(R)\). Moreover, for any fixed center in
\(\mathcal Z\), comparability of \(\ell_N\) with \(r_*(\zeta_N)\) restricts
associated terminals to boundedly many dyadic scales, and net separation bounds
the number of terminals per scale. Bounded overlap of the doubled intrinsic
balls \(B_{2r_*(\zeta)}(\zeta)\) holds within each scale class of
\cite[Definition~5.4]{FRS} by disjointness of the undoubled balls, whereas
doubled balls from distinct classes are mutually disjoint because larger-scale
centers eliminate smaller-scale candidates within four radii. This proves~(i).

\vskip0.13in

\smallskip
\noindent\emph{Step 3: Slab trapping and deep positive caps around \(z\).}
Let \(0<r\le R/32\) with \(r_*(z)/r\) small. By \cite[Lemma~5.6]{FRS}, there
exists \(e\in\mathbb S^{n-1}\) such that \(|t|=|(x-z)\cdot e|\) satisfies
\(\|u-|t|\|_{L^\infty(B_{32r})}\le\eta r\), where \(\eta\to0\) as
\(r_*(z)/r\to0\). Suppose \(\overline N\) intersects \(B_r\). If
\(\ell_N\ge r\), then \(|\zeta_N-z|\le 6\ell_N\), and combining Step~2 with
\cite[Lemma~5.7]{FRS} at \(\zeta_N\) would force \(r_*(z)\ge c_\delta\ell_N\ge c_\delta r\),
impossible once \(r_*(z)/r\) is small. Thus \(\ell_N<r\) and
\(|\zeta_N-z|\le 6r\), whence \cite[Lemma~5.8]{FRS} at scale \(16r\) gives
\(r_*(\zeta_N)\le 2r\). Recentering the lower wedge barrier at \(\zeta_N\) on
\(B_{4r_*(\zeta_N)}(\zeta_N)\subset B_{32r}\) incurs an error at most
\(2\eta r\), so \cite[Lemma~9.8]{FRS} yields \(\ell_N\le C_\delta\eta r\).
Since \(x_N\in B_{2r}\cap\FB(u)\) also satisfies \(|(x_N-z)\cdot e|\le\eta r\),
both deep caps \(B_r\cap\{\pm t>C'_\delta\eta r\}\) lie strictly inside
\(\Omega\) and avoid every closed neck terminal.

Following the local-child nested cover \cite[(9.5)--(9.6)]{FRS}, there is a
descending chain of tree balls keeping \(z\) inside their quarter-balls that
must end at a neck terminal \(N_*\ni z\) (as \(z\in\mathcal Z\) prevents
regular termination). Since \(\ell_{N_*}\le C_\delta\eta r\ll r\), this chain
contains a branching ball \(B\) with \(r/8<r_B\le r/4\) and \(z\in\frac14 B\).
Recentering \(|t|\) at \(x_B\in\FB(u)\) on \(16B\) and invoking affine
stability of optimal wedges \cite[Lemma~9.5]{FRS} aligns \(t\) (after choosing
the sign of \(e\)) with the signed affine representative \(\ell_B\) of \(B\),
so that \(\|\ell_B-t\|_{L^\infty(B)}\le C\eta r\). The two test points
\(x_B\pm\frac12 r_B e\) belong to the respective signed cores of \(B\)
\cite[(9.17)]{FRS} and have signed heights at least \(r/16-\eta r\). Because
each deep cap in \(B_r\) is connected, avoids all neck terminals, and is
covered by \(U_+\cup U_-\), meeting the corresponding core places the entire
cap inside \(U_\pm\). Disjointness of \(U_+\) and \(U_-\) then implies the
complementary upper slab inclusions, proving \eqref{eq:slab}.

\vskip0.13in

\smallskip
\noindent\emph{Step 4: Perimeter and volume estimates on the residual neck set.}
Because the tree is finite and \(U_\pm\) are constructed from \(\Omega\)
(which has smooth boundary) via finitely many spherical and planar cutoffs,
\(U_\pm\) have locally finite perimeter. Inside \(B_{R/8}\), the union of
closed neck terminals \(\overline N\) covers \(Y\) by \cite[(9.20)]{FRS}; any
point of \(\Omega\) outside that closed union possesses an open neighborhood in
\(\Omega\) partitioned into the disjoint open sets \(U_+\) and \(U_-\), and so
cannot belong to \(\partial^*U_\pm\cap\Omega=G_\pm\).

\vskip0.13in

To bound \(\mathcal H^{n-1}(G_\pm\cap\overline N\cap B_{R/8})\), we localize
the face-counting argument of \cite[Lemma~10.3]{FRS}. Suppose a cutoff boundary
face of the signed core of a tree ball \(B\) intersects \(\overline N\). If
\(B\) is a branching ball, any point \(y\) on its artificial cutoff face
satisfies \(|\ell_B(y)|\ge r_B/16\); if \(\ell_N<r_B/4\), evaluating the slab
bound at \(x_N\in\FB(u)\) yields \(r_B/16\le C\delta r_B+\ell_N\), forcing
\(r_B\lesssim_\delta\ell_N\). If instead \(B\) is a regular terminal, then
\(4B\cap\mathcal Z=\varnothing\) while \(|\zeta_N-x_B|\le r_B+5\ell_N\), which
forces \(r_B+5\ell_N\ge 4r_B\) and hence \(\ell_N\ge\frac35 r_B\). Conversely,
when \(r_B\le\ell_N\), let \(A\) be \(B\) (if branching) or its parent (if
regular), and pick \(\zeta_A\in\mathcal Z\cap 4A\); then
\(|\zeta_A-\zeta_N|\le C\ell_N\) and \(r_*(\zeta_A)\le C\delta r_A\), whereas
scale comparison \cite[Lemma~5.7]{FRS} at \(\zeta_N\) gives
\(r_*(\zeta_A)\ge c_\delta\ell_N\), so \(r_B\gtrsim_\delta\ell_N\). Thus
\(r_B\asymp_\delta\ell_N\) for every tree ball contributing a face to
\(\overline N\). Net separation bounds the number of such balls to \(O_\delta(1)\),
and each spherical or affine cutoff face has \((n-1)\)-measure \(O(r_B^{n-1})\),
yielding
\[
  \mathcal H^{n-1}(G_\pm\cap\overline N\cap B_{R/8})\le C_\delta\ell_N^{n-1},
\]
including any portion of \(G_\pm\) lying on \(\partial N\). Summing this bound
and \(|\overline N|=C_n\ell_N^n\) over all terminals meeting \(B_r\), and
factoring out one or two powers of \(\ell(r)\), gives \eqref{eq:area-volume}.$\square$
\end{proof}

\vskip0.33in

\subsection{Sternberg--Zumbrun projection and neck content control}\label{sec:sz}

To control the \((n-2)\)-dimensional neck content \(Q(s)\) in terms of how much
\(\nabla u\) deviates from a two-dimensional plane, we project the translation
Jacobi fields onto an \((n-2)\)-dimensional subspace \(M\subset\R^n\) and apply
a Sternberg--Zumbrun-type stability estimate. The key geometric input is that
vanishing of the projected directional curvature at a neck scale splits off an
\(\R^{n-3}\) Euclidean factor, reducing the limit to a three-dimensional stable
solution ruled out by \cite{CFFS}.

\vskip0.13in

\begin{lemma}\label{lem:gap}
Let \(M\subset\R^n\) be an \((n-2)\)-dimensional linear subspace, denote by
\(P_M\) the orthogonal projection onto \(M\), put \(F=P_M\nabla u\), and define
the projected directional curvature density
\[
  \mathcal A_M^2=|F|^2\bigl|D(F/|F|)\bigr|^2\quad\text{where }F\ne0,
  \qquad \mathcal A_M^2=0\quad\text{where }F=0.
\]
Then there exists a constant \(c=c(n,\eta_0)>0\), independent of \(M\), such
that for every intrinsic center \(\zeta\in\mathcal Z\),
\begin{equation}\label{eq:gap}
  \int_{B_{2r_*(\zeta)}(\zeta)\cap\Omega}\mathcal A_M^2\dd x
  \ge c\,r_*(\zeta)^{n-2}.
\end{equation}
Moreover, for every smooth cutoff \(\eta\in C_c^\infty(\R^n)\),
\begin{equation}\label{eq:SZ}
  \int_\Omega\mathcal A_M^2\eta^2\dd x
  \le \int_\Omega |P_M\nabla u|^2|\nabla\eta|^2\dd x.
\end{equation}
Consequently, in the setting of Lemma~\ref{lem:fixed-root-tree}, whenever
\(0<s\le R/32\) and \(s<r\le R\) are such that all doubled intrinsic balls
\(B_{2r_*(\zeta_N)}(\zeta_N)\) associated with terminals whose closures meet
\(B_s\) lie inside \(B_{r/4}\), one has
\begin{equation}\label{eq:projected-charge}
  Q(s)\le C r^{-2}\int_{B_{r/2}\cap\Omega}|P_M\nabla u|^2\dd x.
\end{equation}
\end{lemma}

\begin{proof}
With respect to any orthonormal basis of \(M\), the \(n-2\) components of
\(F=P_M\nabla u\) are directional derivatives of \(u\), hence classical Jacobi
fields annihilated by the linearized operator inside \(\Omega\) and satisfying
the linearized Robin condition on \(\FB(u)\). Applying the vector-valued
Sternberg--Zumbrun stability inequality for finite families of Jacobi fields
\cite[Lemma~6.3]{FRS} directly gives \eqref{eq:SZ}.

\vskip0.13in

To establish the uniform lower bound \eqref{eq:gap}, suppose by contradiction
that along a sequence of normalized solutions \(u_k\), neck centers
\(\zeta_k\), and \((n-2)\)-planes \(M_k\), the rescaled integral
\(r_*(\zeta_k)^{2-n}\int_{B_{2r_*(\zeta_k)}(\zeta_k)\cap\Omega_k}\mathcal A_{M_k}^2\)
tends to zero. Set \(\tilde u_k(x)=r_*(\zeta_k)^{-1}u_k(\zeta_k+r_*(\zeta_k)x)\).
By definition of \(r_*(\zeta_k)\), each \(\tilde u_k\) has \(L^n\)-Hessian
mass \(\eta_0^n\) in \(B_1\), while scale comparison \cite[Lemma~5.7]{FRS}
provides uniform Hessian bounds on every fixed ball. Extracting a diagonal
subsequence via neck compactness \cite[Lemmas~4.6 and~6.8]{FRS} and
compactness of the Grassmannian, we obtain an \((n-2)\)-plane \(M\) and an
entire classical stable limit \(v\) satisfying
\[
  \int_{B_1\cap\{v>0\}}|D^2v|^n\dd x=\eta_0^n>0.
\]
On the open subset of \(B_2\cap\{v>0\}\) where \(F=P_M\nabla v\ne0\), smooth
interior convergence implies pointwise convergence of
\(\mathcal A_{M_k}^2(\tilde u_k)\) to
\(\mathcal A_M^2(v)=\bigl|(I-\widehat F\otimes\widehat F)DF\bigr|^2\), where
\(\widehat F=F/|F|\). Fatou's lemma gives vanishing almost everywhere on \(\{F\ne0\}\cap B_2\cap\{v>0\}\).
Continuity on this set and the convention \(\mathcal A_M^2=0\) on \(\{F=0\}\)
then give pointwise vanishing on \(B_2\cap\{v>0\}\).

\vskip0.13in

Now let \(D\) be any connected component of \(\{v>0\}\) intersecting \(B_1\).
If \(F\) is not identically zero on \(D\cap B_2\), pick a small interior ball
\(B'\Subset D\cap B_2\) on which \(F\ne0\); vanishing of \(\mathcal A_M^2(v)\)
implies that \(\widehat F\) is constant on \(B'\), so all directional
derivatives of \(v\) along the \((n-3)\)-dimensional subspace
\(W=M\cap\widehat F^\perp\) vanish on \(B'\), and hence throughout \(D\) by
unique continuation for harmonic functions. If instead \(F\equiv0\) on
\(D\cap B_2\), the same unique continuation holds for any \((n-3)\)-dimensional
subspace \(W\subset M\). Because \(v\) is constant and strictly positive along
every line parallel to \(W\) passing through a point of \(D\), no such line can
meet \(\partial D\). In orthogonal coordinates adapted to \(\R^n=\R^3\times W\),
the component splits metrically as a cylinder \(D=D_0\times\R^{n-3}\) with
\(v(x,y)=v_0(x)\). Since the fibers \(\{x\}\times\R^{n-3}\) are tangent to
\(\partial D\), smooth local one-sidedness descends to \(D_0\subset\R^3\), so
extending \(v_0\) by zero outside \(D_0\) produces an entire classical
solution to \eqref{eq:bernoulli} in \(\R^3\).

\vskip0.13in

Furthermore, \(v|_D\) inherits stability on its own closure: every point of
\(\overline D\) is either interior to \(D\) or admits a smooth one-sided
free-boundary chart disjoint from all other components of \(\{v>0\}\), allowing
any compactly supported test function on \(\overline D\) to be extended by zero
to the remaining components without changing \(\mathcal Q_v\). Given an
admissible test function \(\phi\) on \(\overline{D_0}\subset\R^3\) and a
nontrivial cutoff \(\eta\in C_c^\infty(\R^{n-3})\), testing with the product
\(\Phi_A(x,y)=\phi(x)\eta(y/A)\) yields, since the Euclidean factor carries
zero principal curvatures,
\[
  0\le \frac{\mathcal Q_{v|_D}[\Phi_A;\R^n]}{\|\eta(\cdot/A)\|_{L^2}^2}
  =\mathcal Q_{v_0}[\phi;\R^3]
   +A^{-2}\frac{\|\nabla\eta\|_{L^2}^2}{\|\eta\|_{L^2}^2}\int_{D_0}\phi^2\dd x.
\]
Letting \(A\to\infty\) shows that \(v_0\) is a classical stable solution in
\(\R^3\). By three-dimensional stable rigidity \cite[Theorem~1.5]{CFFS},
\(D^2v_0\equiv0\) on \(D_0\), and therefore \(D^2v\equiv0\) on every component
of \(\{v>0\}\) meeting \(B_1\), contradicting
\(\int_{B_1\cap\{v>0\}}|D^2v|^n\dd x=\eta_0^n>0\). This proves \eqref{eq:gap}.

\vskip0.13in

Finally, under the containment hypothesis of \eqref{eq:projected-charge}, we
sum \eqref{eq:gap} over the distinct centers \(\zeta_N\) associated with
terminals whose closures meet \(B_s\). Using \(\ell_N\asymp_\delta r_*(\zeta_N)\), the
bounded multiplicity of \(N\mapsto\zeta_N\), and the bounded overlap of
\(B_{2r_*(\zeta_N)}(\zeta_N)\) from Lemma~\ref{lem:fixed-root-tree}(i), we
bound \(Q(s)\) by a constant multiple of
\(\int_{B_{r/4}\cap\Omega}\mathcal A_M^2\dd x\), and applying \eqref{eq:SZ}
with a smooth cutoff \(\eta\equiv1\) on \(B_{r/4}\) supported in \(B_{r/2}\)
with \(|\nabla\eta|\le C/r\) yields \eqref{eq:projected-charge}.$\square$
\end{proof}

\vskip0.33in

\subsection{Fixed-center excess decay and proof of Theorem~\ref{thm:transfer}}\label{sec:contraction}

Retaining the localized decomposition of Lemma~\ref{lem:fixed-root-tree}
around \(z\in\mathcal Z\), for any radius \(0<r\le R/32\) and any pair
\(L=(L_+,L_-)\) of affine functions \(L_\pm(x)=b_\pm\cdot(x-z)+d_\pm\) with
unit slopes \(|b_\pm|=1\), we introduce the coupled quadratic-neck excess
\begin{equation}\label{eq:excess}
  \Exc(r,L)^2
  =r^{-n-2}\sum_\pm\int_{U_\pm\cap B_r}(u-L_\pm)^2\dd x
   +\kappa\,r^{2-n}Q(r),
\end{equation}
where \(\kappa>0\) is a dimensional coupling weight to be fixed below. The
heart of the rigidity proof is the following excess-decay contraction.

\vskip0.33in

\begin{lemma}\label{lem:contract}
There exist constants \(\theta\in(0,1/32)\), \(\kappa>0\), and \(\lambda_*>0\),
depending only on \(n\) and the fixed tree aperture \(\delta\), such that
whenever \(0<r\le R/32\) and \(\Exc(r,L)\le\lambda_*\), there exists a new pair
\(L'=(L'_+,L'_-)\) of unit-slope affine functions satisfying
\[
  \Exc(\theta r,L')\le \tfrac12\Exc(r,L).
\]
\end{lemma}

\begin{proof}
To make it transparent that no circularity occurs among the constants, we note
at the outset that \(\delta\) has already been fixed in
Lemma~\ref{lem:fixed-root-tree}; below we will first select \(\theta\in(0,1/32)\)
solely from interior harmonic Taylor decay, next fix \(\kappa>0\) small
relative to \(\theta^{n-2}\), and finally choose \(\lambda_*>0\) small enough
for all compactness and higher-order absorption steps at the fixed pair
\((\theta,\kappa)\).

\vskip0.13in

\smallskip
\noindent\emph{Step 1: Rescaling, neck smallness, and approximate normal alignment.}
Translating \(z\) to the origin and rescaling from \(B_r(z)\) to \(B_1(0)\),
put \(\lambda=\Exc(1,L)\le\lambda_*\), and let \(Q,S,\ell\) refer to the
rescaled tree. Directly from \eqref{eq:excess} and \eqref{eq:area-volume},
setting \(q=\lambda^2/\kappa\), we have
\begin{equation}\label{eq:small-A}
  Q(1)\le q,\qquad
  \ell(1)\le q^{1/(n-2)},\qquad
  S(1)\le C q^{1+1/(n-2)},\qquad
  |Y\cap B_1|\le C q^{1+2/(n-2)}.
\end{equation}
For each fixed \(\kappa>0\), these bounds give
\[
  \frac{S(1)+|Y\cap B_1|}{\lambda^2}
  \le \frac{C}{\kappa}\bigl(q^{1/(n-2)}+q^{2/(n-2)}\bigr)\longrightarrow0
  \qquad (\lambda\to0).
\]
We may therefore impose \(S(1)+|Y\cap B_1|\le\lambda^2\) by shrinking
the smallness threshold after fixing \(\kappa\). Moreover, the terminal \(N_*\) containing
the origin satisfies \(|\zeta_{N_*}|\le 5\ell_{N_*}\) and
\(r_*(\zeta_{N_*})\asymp\ell_{N_*}\le\ell(1)\), so scale comparison
\cite[Lemma~5.7]{FRS} gives \(r_*(0)\le C\ell(1)\to0\) in the rescaled
variables. Lemma~\ref{lem:fixed-root-tree}(ii) therefore applies on \(B_1\);
rotating coordinates so that \(e=e_n\), the slab inclusions \eqref{eq:slab}
hold with \(t=x_n\), and \(L^2\) comparison of \(L_\pm(x)=b_\pm\cdot x+d_\pm\)
against \(\pm x_n\) on fixed interior balls inside the deep caps gives
\(b_\pm=\pm e_n+o(1)\) and \(d_\pm=o(1)\) as \(\lambda\to0\).

\vskip0.13in

\smallskip
\noindent\emph{Step 2: One-sided lower barrier and Caccioppoli energy bound.}
Fix a sign \(\sigma\in\{+,-\}\) and write \(w=u-L_\sigma\) and \(h=\nabla w\)
on \(U_\sigma\). We first show that \(w\ge -C\lambda\) throughout
\(B_{2/3}\) (with \(w\) defined globally via \(u\ge0\)). Indeed, at any point
\(x\in B_{2/3}\) where \(w(x)<0\), necessarily \(L_\sigma(x)>0\), which by
Step~1 forces \(\sigma x_n\ge -o(1)\). Shifting along the unit slope to
\(x'=x+\frac14 b_\sigma\), the global \(1\)-Lipschitz bound \(|\nabla u|\le1\)
implies \(u(x)-u(x')\ge -|x-x'|=-\frac14=L_\sigma(x)-L_\sigma(x')\), so
\(w(x)\ge w(x')\). Because \(|x'|\le\frac{2}{3}+\frac{1}{4}=\frac{11}{12}\) and
\(\sigma x'_n\ge\frac14-o(1)\), a ball of fixed radius around \(x'\) is
contained in the deep cap \(U_\sigma\cap B_1\) by \eqref{eq:slab}, where \(w\)
is harmonic and the mean-value inequality yields \(|w(x')|\le C\lambda\). Thus
\(w+C\lambda\ge0\) on \(B_{2/3}\).

\vskip0.13in

Along the physical free boundary \(\partial^*U_\sigma\cap\FB(u)\), the outward
unit normal to \(U_\sigma\subset\Omega\) coincides \(\mathcal H^{n-1}\)-almost
everywhere with \(\nu=-\nabla u\), and since \(|b_\sigma|=|\nabla u|=1\), the
outward normal derivative satisfies
\(\partial_\nu w=-1+b_\sigma\cdot\nabla u\le0\). Choose a smooth cutoff
\(\zeta\in C_c^\infty(B_{2/3})\) with \(\zeta\equiv1\) on \(B_{3/5}\). Because
\(u\) is classical, its positive-phase derivatives extend smoothly across
\(\FB(u)\) on \(\operatorname{supp}\zeta\), permitting application of the
Gauss--Green--De Giorgi formula to \(\zeta^2(w+C\lambda)\nabla w\) on the
Caccioppoli set \(U_\sigma\). The boundary integral along
\(\partial^*U_\sigma\cap\FB(u)\) is nonpositive since \(w+C\lambda\ge0\) and
\(\partial_\nu w\le0\), while along the interior partition boundary \(G_\sigma\)
we have \(|\nabla w|\le 2\) and \(|w+C\lambda|\le C\), contributing at most
\(CS(1)\). Absorbing the interior \(\nabla\zeta\) cross-term via Young's
inequality and invoking \eqref{eq:small-A} gives the Caccioppoli bound
\begin{equation}\label{eq:energy}
  \int_{U_\sigma\cap B_{3/5}}|h|^2\dd x
  \le C\bigl(\lambda^2+S(1)\bigr)\le C\lambda^2,
\end{equation}
The constants in this bound depend only on \(n\) and the fixed aperture
\(\delta\). The smallness threshold may depend on \(\kappa\).

\vskip0.13in

\smallskip
\noindent\emph{Step 3: Pohozaev--Rellich flux identity and \(BV\) trace control.}
Still fixing a side \(\sigma\) with \(b=b_\sigma\), consider the Rellich
stress-energy vector field
\[
  K_b(h)=\tfrac12|h|^2 b-(b\cdot h)h,
\]
which satisfies \(\diver K_b(h)=0\) in \(U_\sigma\) because \(w\) is harmonic
and \(D^2w\) is symmetric. On \(\FB(u)\), writing \(g=\nabla u\) and setting
\(d=1-b\cdot g=1+b\cdot\nu\ge0\), the unit-length identities \(|g|=|b|=1\) give
\(\frac12|h|^2=d\) and \(b\cdot h=-d\), so \(K_b(h)=dg\) and hence
\begin{equation}\label{eq:flux}
  K_b(h)\cdot\nu=\partial_\nu w=-d\le0\quad\text{on }\FB(u).
\end{equation}
Testing \(\diver K_b(h)=0\) on \(U_\sigma\) against a nonnegative cutoff equal
to \(1\) on \(B_{1/2}\) and supported in \(B_{3/5}\), and bounding the interior
gradient term by \eqref{eq:energy} and the co-normal flux across \(G_\sigma\)
by \(CS(1)\), we obtain the boundary alignment estimate
\begin{equation}\label{eq:trace}
  \int_{\partial^*U_\sigma\cap\FB(u)\cap B_{1/2}} d\dd\mathcal H^{n-1}
  \le C\bigl(\lambda^2+S(1)\bigr)\le C\lambda^2.
\end{equation}
Applying Gauss--Green next to \(\phi\nabla w\) for arbitrary
\(\phi\in C_c^\infty(B_{1/2})\), and using \(|\partial_\nu w|=d\) on
\(\partial^*U_\sigma\cap\FB(u)\) together with \(|h|\le2\) on \(G_\sigma\),
\eqref{eq:trace} yields the weak homogeneous Neumann bound
\begin{equation}\label{eq:weak}
  \left|\int_{U_\sigma}\nabla(w/\lambda)\cdot\nabla\phi\dd x\right|
  \le C\lambda\|\phi\|_{L^\infty}
  \qquad\forall\,\phi\in C_c^\infty(B_{1/2}).
\end{equation}

\vskip0.13in

Crucially, \eqref{eq:trace} also prevents \(L^2\) energy from concentrating
inside the vanishing slab \(\{|x_n|\le\delta_1\}\). Pick \(\chi\in C_c^\infty(B_{1/2})\)
with \(0\le\chi\le1\) and \(\chi\equiv1\) on \(B_{3/8}\), and set
\(f=\chi^2 w^2\). By Cauchy--Schwarz, \eqref{eq:excess}, and \eqref{eq:energy},
\(\int_{U_\sigma}|\nabla f|\dd x\le C\lambda^2\). On
\(\partial^*U_\sigma\cap\FB(u)\), decompose \(1=-b\cdot\nu+d\): applying
Gauss--Green to the vector field \(fb\) bounds the integral of \(f(-b\cdot\nu)\)
by \(\int_{U_\sigma}|\nabla f|\dd x+CS(1)\), while boundedness of \(f\) and
\eqref{eq:trace} bound the integral of \(fd\) by \(C\lambda^2\), and the trace
of \(f\) on \(G_\sigma\) contributes at most \(CS(1)\). Consequently, the zero
extension \(f\one_{U_\sigma}\) belongs to \(BV(\R^n)\) with total variation
\(O(\lambda^2)\), and the Euclidean \(BV\)--Sobolev embedding yields the
improved integrability bound
\begin{equation}\label{eq:Lq}
  \|w\|_{L^{2n/(n-1)}(U_\sigma\cap B_{3/8})}\le C\lambda.
\end{equation}

\vskip0.13in

\smallskip
\noindent\emph{Step 4: Harmonic linearization via Neumann reflection.}
We now establish the affine excess improvement by compactness as \(\lambda\to0\).
Along any sequence with \(\lambda\to0\), the slab width in \eqref{eq:slab}
shrinks to zero, so \(U_\pm\cap B_1\) converge to the upper and lower
half-balls \(B_1^\pm=B_1\cap\{\pm x_n>0\}\). For the normalized defects
\(v_\pm=(u-L_\pm)/\lambda\), Hölder's inequality and \eqref{eq:Lq} show that on
any thin neighborhood \(T_\varepsilon=B_{3/8}\cap\{|x_n|<\varepsilon\}\),
\[
  \int_{T_\varepsilon\cap U_\pm}|v_\pm|^2\dd x
  \le C|T_\varepsilon|^{1/n}=O(\varepsilon^{1/n}),
\]
uniformly in \(\lambda\). Away from \(\{x_n=0\}\), standard interior harmonic
estimates give compact convergence on compact subsets of \(B_{3/8}^\pm\), so
letting \(\varepsilon\downarrow0\) upgrades convergence of the zero extensions
\(v_\pm\one_{U_\pm}\) to strong convergence in \(L^2(B_{1/3})\) toward harmonic
limits \(\bar v_\pm\) on \(B_{1/3}^\pm\) satisfying
\(\|\bar v_\pm\|_{L^2(B_{1/3}^\pm)}\le1\). Simultaneously, \eqref{eq:energy}
provides weak \(L^2\) convergence of \((\nabla v_\pm)\one_{U_\pm}\) to
\((\nabla\bar v_\pm)\one_{B_{1/3}^\pm}\), and passing to the limit in
\eqref{eq:weak} shows that each \(\bar v_\pm\) satisfies the weak homogeneous
Neumann boundary condition \(\partial_n\bar v_\pm=0\) along
\(B_{1/3}\cap\{x_n=0\}\). Even reflection across \(\{x_n=0\}\) therefore
extends \(\bar v_\pm\) harmonically across the nodal hyperplane to all of
\(B_{1/3}\).

\vskip0.13in

Let \(p_\pm(x)\) be the first-order Taylor polynomials of these even-reflected
harmonic functions at the origin; by evenness in \(x_n\), \(\partial_n p_\pm=0\),
so \(\nabla p_\pm\) are purely tangential and bounded by a dimensional
constant, while interior harmonic Taylor estimates give
\(\sum_\pm\|\bar v_\pm-p_\pm\|_{L^2(B_\theta^\pm)}^2\le C_0\theta^{n+4}\) for
all \(\theta\in(0,1/32)\). Define updated unit-slope affine maps \(L'_\pm\) by
normalizing the slope of \(L_\pm+\lambda p_\pm\) while keeping its value at the
origin: because \(b_\pm\to\pm e_n\) and \(\partial_n p_\pm=0\), we have
\[
  \bigl|b_\pm+\lambda\nabla p_\pm\bigr|-1
  =\lambda\,b_\pm\cdot\nabla p_\pm+O(\lambda^2)=o(\lambda),
\]
so slope normalization perturbs \(L_\pm+\lambda p_\pm\) on \(B_\theta\) by at
most \(o(\lambda)\theta\). Fix \(\theta\in(0,1/32)\) small enough that
\(C_0\theta^2<1/16\). By strong \(L^2\) convergence (via a contradiction
sequence if the bound failed for arbitrarily small \(\lambda\)), for all
sufficiently small \(\lambda\) we obtain
\begin{equation}\label{eq:affine-next}
  \theta^{-n-2}\sum_\pm\int_{U_\pm\cap B_\theta}(u-L'_\pm)^2\dd x
  \le \tfrac18\lambda^2.
\end{equation}

\vskip0.13in

\smallskip
\noindent\emph{Step 5: Controlling the terminal neck content and closing the contraction.}
Choose an \((n-2)\)-dimensional subspace \(M\subset b_+^\perp\cap b_-^\perp\).
On \(U_\pm\) we have \(P_M\nabla u=P_M(\nabla u-b_\pm)\), so
\(|P_M\nabla u|\le|h|\), whereas on the residual set \(Y\) we simply use
\(|P_M\nabla u|\le|\nabla u|\le1\). Combining \eqref{eq:energy} with the
volume bound \(|Y\cap B_1|=o(\lambda^2)\) from \eqref{eq:small-A} yields
\(\int_{B_{1/2}\cap\Omega}|P_M\nabla u|^2\dd x\le C\lambda^2\). Moreover, for
small \(\lambda\), every neck terminal meeting \(B_\theta\) has doubled
intrinsic ball contained in \(B_{\theta+C\ell(1)}\subset B_{1/4}\). Applying
\eqref{eq:projected-charge} from Lemma~\ref{lem:gap} (at \(s=\theta\) and
\(r=1\) in the rescaled variables) gives
\(Q(\theta)\le C_1\lambda^2\), with \(C_1\) depending only on \(n\) and the
fixed aperture \(\delta\).

\vskip0.13in

We now fix \(\kappa>0\) so that \(\kappa C_1\theta^{2-n}\le 1/8\), and then fix
\(\lambda_*>0\) small enough for all smallness steps above at this fixed
\((\theta,\kappa)\). Adding \(\kappa\theta^{2-n}Q(\theta)\le\frac18\lambda^2\)
to \eqref{eq:affine-next} gives
\(\Exc(\theta,L')^2\le\frac14\lambda^2=\frac14\Exc(1,L)^2\), completing the
proof. $ \square$
\end{proof}

\newpage

\begin{proof}[Proof of Theorem~\ref{thm:transfer}]
Fix \(z\in\mathcal Z\), build the localized tree of
Lemma~\ref{lem:fixed-root-tree} at root radius \(R\), and set \(r_0=R/32\). By
Lemma~\ref{lem:fixed-root-tree}(i), \(\max_N\ell_N=o(R)\), so for all large
\(R\) the doubled intrinsic balls associated with terminals meeting \(B_{r_0}\)
lie inside \(B_{R/4}\). Let \(e_R\in\mathbb S^{n-1}\) be an approximate wedge
direction at scale \(R\) given by \cite[Lemma~5.6]{FRS}, and choose any
\((n-2)\)-plane \(M_R\subset e_R^\perp\). Rescaling \(B_R(z)\) to \(B_1(0)\)
and rotating \(e_R\) to \(e_n\), the rescaled functions converge uniformly to
\(|x_n|\) on \(B_1\), so interior harmonic derivative convergence off
\(\{x_n=0\}\) combined with \(|\nabla u|\le1\) gives
\(R^{-n}\int_{B_{R/2}\cap\Omega}|P_{M_R}\nabla u|^2\dd x=o(1)\) as
\(R\to\infty\). Applying \eqref{eq:projected-charge} with \(s=r_0\) and \(r=R\)
therefore yields the initial neck bound
\begin{equation}\label{eq:seed}
  r_0^{2-n}Q(r_0)=32^{n-2}R^{2-n}Q(r_0)=o(1)\qquad (R\to\infty).
\end{equation}
Take \(L_{0,\pm}=\pm t\) for the wedge coordinate \(t=(x-z)\cdot e\) from
\eqref{eq:slab} at radius \(r_0\). Since \(\pm t> -\delta_{r_0}r_0\) on
\(U_\pm\cap B_{r_0}\), we have
\[
  \sup_{U_\pm\cap B_{r_0}}|u-L_{0,\pm}|
  \le \|u-|t|\|_{L^\infty(B_{r_0})}+2\delta_{r_0}r_0=o(r_0).
\]
The normalized affine term in \eqref{eq:excess} thus tends to zero.
Together with \eqref{eq:seed}, this gives \(\Exc(r_0,L_0)\le\lambda_*\)
once we fix a sufficiently large root radius \(R\).

\vskip0.13in

Iterating Lemma~\ref{lem:contract} on this fixed tree along the geometric
scales \(r_j=\theta^j r_0\) produces unit-slope pairs \(L_j\) satisfying
\(\Exc(r_j,L_j)\le 2^{-j}\Exc(r_0,L_0)\). Because the fixed neck terminal
\(N_*\) containing \(z\) has a strictly positive radius \(\ell_*>0\) and meets
\(B_{r_j}\) for every \(j\ge0\), we arrive at
\[
  0<\ell_*^{n-2}\le Q(r_j)
  \le \kappa^{-1}r_j^{n-2}\Exc(r_j,L_j)^2
  \le \kappa^{-1}r_j^{n-2}4^{-j}\Exc(r_0,L_0)^2\longrightarrow0
  \quad (j\to\infty),
\]
a contradiction. Thus no nonflat entire classical stable solution exists, so
every entire classical stable solution satisfies \(D^2u\equiv0\) in \(\Omega\).

\vskip0.13in

Finally, when \(D^2u\equiv0\) in \(\Omega\), \(\nabla u\) is constant on each
connected component of \(\Omega\). If \(\nabla u\equiv0\) on a component, that
component has empty boundary inside \(\R^n\), so \(\Omega=\R^n\) and \(u\) is a
positive constant. Otherwise, each component has constant nonzero gradient---of
unit length whenever \(\FB(u)\ne\varnothing\)---and since \(u=0\) on its
finite boundary, the component must be an open affine half-space. At most two
disjoint open half-spaces can exist in \(\R^n\), and if there are two, their
normals must be opposite. Thus \(u\) is one-dimensional in every case.
$ \square$\end{proof}

\vskip0.33in

\subsection{Classification up to dimension six and sharpness in dimension seven}\label{subsec:sharp-dim}

\begin{proof}[Proof of Corollary~\ref{cor:classification}]
By the smooth stable cone classification \cite[Theorem~3 and Remark~4]{survey},
which incorporates \cite{OAI} in dimensions five and six, and the minimizing
cone of De Silva and Jerison \cite{DSJ}, one has
\(n_{\mathrm{cone}}^*>6\). Theorem~\ref{thm:transfer} therefore yields
\(D^2u\equiv0\) and one-dimensionality in dimensions \(4\le n\le6\), while
dimensions \(n=2,3\) are established in \cite{CFFS,survey}. The affine
component analysis at the end of the proof of Theorem~\ref{thm:transfer},
together with \(|\nabla u|=1\) on \(\FB(u)\) and smooth local one-sidedness of
\(\Omega\) (which precludes contact of two parallel free-boundary hyperplanes
and forces \(b<a\)), yields the stated family of profiles. Conversely, every
listed profile has flat free boundary with \(H\equiv0\), so
\(\mathcal Q_u[\phi;\R^n]=\int_\Omega|\nabla\phi|^2\dd x\ge0\) and \(u\) is
stable. $ \square$
\end{proof}

\vskip0.33in

\begin{corollary}\label{cor:critical}
The smallest dimension admitting a nonflat entire classical stable solution to
\eqref{eq:bernoulli} is \(n_{\mathrm{entire}}^*=7\).
\end{corollary}

\begin{proof}
Corollary~\ref{cor:classification} rules out dimensions \(2\le n\le6\), and in
dimension \(n=1\) harmonicity trivially gives \(u''=0\) on every positive
interval. In dimension \(n=7\), let \(C_0\) be a nonflat one-homogeneous global
minimizer constructed by De Silva and Jerison \cite{DSJ}. Applying
the Hardt--Simon-type foliation theorem of Edelen, Spolaor and Velichkov
\cite[Theorem~1.5]{ESV} to \(C_0\) yields an entire global minimizer \(V\) in
\(\R^7\) with empty singular set whose blow-down at infinity is \(C_0\).
Smoothness of \(\FB(V)\) and global minimality imply that \(V\) is an entire
classical stable solution. Since the affine-component classification of
solutions with \(D^2u\equiv0\) in the proof of Theorem~\ref{thm:transfer} holds
in arbitrary dimension, any solution with vanishing positive-phase Hessian has
blow-down equal to \(0\), a half-space profile, or a unit wedge; because \(C_0\)
is none of these, \(V\) is genuinely nonflat. $ \square$
\end{proof}

\vskip0.33in

\section{\bf Local curvature estimates and compactness of classical solutions}\label{sec:curvature}

In this section we assemble the local regularity, multi-sheet graphical
convergence, and point-picking tools needed both to deduce the interior
curvature bound (Corollary~\ref{cor:curv}) and to extract finite-Morse-index
bubbles in Section~\ref{sec:index-concentration}. Throughout what follows, a
function of the form \(|(x-p)\cdot e|\) with \(e\in\mathbb S^{n-1}\) is called
a \emph{unit wedge} (and a \emph{centered unit wedge} when \(p=0\)); note that
its classical Hessian vanishes identically on its positive phase
\(\R^n\setminus\{x:(x-p)\cdot e=0\}\).

\vskip0.13in

\subsection{Interior Lipschitz bounds and graphical sheet convergence}\label{subsec:sheets}

We begin with a universal interior Lipschitz estimate for classical solutions
whose free boundary meets an inner ball, independent of stability or index.

\begin{lemma}\label{lem:local-lipschitz}
Let \(u\) be a classical solution to \eqref{eq:bernoulli} in \(B_R\subset\R^n\)
such that \(\FB(u)\cap B_{R/2}\ne\varnothing\). For every \(\theta\in(0,1)\),
\[
  R^{-1}\|u\|_{L^\infty(B_{\theta R})}
  +\|\nabla u\|_{L^\infty(B_{\theta R}\cap\{u>0\})}
  \le C(n,\theta).
\]
\end{lemma}

\begin{proof}
Rescale to \(R=1\) and write \(d(x)=\operatorname{dist}(x,\{u=0\})\). Fix
\(\theta\in[1/2,1)\) and set \(\delta_0=\frac14(1-\theta)\). At any point
\(x\in B_\theta\cap\{u>0\}\) with \(d(x)\le\delta_0\), the ball \(B_{d(x)}(x)\)
lies in \(B_1\cap\{u>0\}\) and touches \(\FB(u)\) at some contact point \(q\).
By Harnack's inequality, \(u\ge c_n u(x)\) on \(\partial B_{d(x)/2}(x)\);
comparing \(u\) on the annulus \(B_{d(x)}(x)\setminus\overline{B_{d(x)/2}(x)}\)
with the radial harmonic barrier taking values \(c_n u(x)\) on the inner sphere
and \(0\) on the outer sphere, the Hopf boundary condition \(|\nabla u(q)|=1\)
yields \(1\ge c'_n u(x)/d(x)\), and interior harmonic gradient estimates on
\(B_{d(x)/2}(x)\) then give \(|\nabla u(x)|\le C_n\). When \(d(x)>\delta_0\),
join a point of \(\FB(u)\cap B_{1/2}\) to \(x\) by a line segment: after its
last exit from \(\{d\le\delta_0\}\), the remaining segment stays at distance at
least \(\delta_0\) from \(\{u=0\}\cup\partial B_1\) and is covered by a Harnack
chain of length \(O_\theta(1)\), bounding \(u(x)\) and hence \(|\nabla u(x)|\).
$\square$
\end{proof}

\vskip0.23in

\begin{remark}\label{rem:inward-segment}
Suppose \(u\) is a classical solution whose positive-phase Hessian satisfies
\(|D^2u|\le K\) in a neighborhood of \(q\in\FB(u)\), and let
\(\nu_{\mathrm{in}}=\nabla u(q)\) be the inward unit normal to \(\Omega\) at
\(q\). Integrating along \(t\mapsto q+t\nu_{\mathrm{in}}\) up to any putative
first exit from \(\Omega\) gives, for all \(0<t\le 1/K\) within the
neighborhood,
\begin{equation}\label{eq:inward-segment}
  u(q+t\nu_{\mathrm{in}})\ge t-\tfrac12 K t^2\ge \tfrac12 t,
\end{equation}
which rules out any exit on that interval and guarantees an inward positive
segment of uniform length. Furthermore, for any locally uniformly bounded
family of classical solutions sharing a uniform positive-phase Hessian bound on
an open set, \(|\nabla u|\) is locally uniformly bounded as well: near
\(\FB(u)\) one integrates \(D^2u\) from a nearest free-boundary point (where
\(|\nabla u|=1\)), away from \(\{u=0\}\) one applies interior harmonic
estimates, and integrating along positive subintervals of arbitrary line
segments extends this gradient bound to a uniform ambient Lipschitz constant.
\end{remark}

\vskip0.23in

\begin{lemma}\label{lem:hessian-convergence}
Suppose classical solutions \(v_j\) converge locally uniformly on an open set
\(G\subset\R^n\) to a function \(v\), and suppose their positive-phase Hessians
\(D^2v_j\) are uniformly bounded on every compact subset of \(G\). Then for
every \(1\le p<\infty\),
\[
  D^2v_j\one_{\{v_j>0\}}\longrightarrow D^2v\one_{\{v>0\}}
  \quad\text{strongly in }L^p_{\mathrm{loc}}(G).
\]
\end{lemma}

\begin{proof}
Let \(F\subset G\) be the set of subsequential accumulation points in \(G\) of
the free boundaries \(\FB(v_j)\). Clearly \(F\) is relatively closed in \(G\)
and contained in \(\{v=0\}\). Given \(q\in F\), pass to a subsequence with
\(q_j\in\FB(v_j)\) converging to \(q\) and inward unit normals converging to
\(\nu_{\mathrm{in}}\). For every sufficiently small \(r>0\), passing to the
limit in \eqref{eq:inward-segment} at \(t=r/2\) gives
\(v(q+\frac12 r\nu_{\mathrm{in}})\ge r/4\); combined with the uniform local
Lipschitz bound from Remark~\ref{rem:inward-segment}, this produces a ball of
radius \(cr\) inside \(B_r(q)\cap\{v>0\}\). Hence \(F\) is locally porous and
has Lebesgue measure zero. On any small ball \(B'\Subset G\setminus F\),
\(\FB(v_j)\cap B'\) is empty for all large \(j\), so each \(v_j\) is harmonic
on \(B'\) (possibly identically zero), and interior harmonic convergence shows
that \(D^2v_j\one_{\{v_j>0\}}\to D^2v\one_{\{v>0\}}\) uniformly on compact
subsets of \(G\setminus F\) (using the minimum principle to distinguish
\(v>0\) from \(v\equiv0\) on \(B'\)). Since \(|F|=0\), dominated convergence
yields strong convergence in \(L^p_{\mathrm{loc}}(G)\).$\square$
\end{proof}

\vskip0.23in

\begin{lemma}\label{lem:entire-slope}
Let \(v_j\) be classical solutions on expanding balls \(B_{R_j}\) (\(R_j\to\infty\))
with a uniform Lipschitz bound \(\|\nabla v_j\|_{L^\infty(B_{R_j})}\le C_0\). If
\(v_j\to v\) locally uniformly on \(\R^n\), then \(|\nabla v|\le1\) throughout
\(\{v>0\}\).
\end{lemma}

\begin{proof}
By interior harmonic convergence, \(v\) is harmonic on \(\{v>0\}\). If
\(\{v>0\}=\R^n\), Liouville's theorem for nonnegative harmonic functions
implies that \(v\) is constant and \(\nabla v\equiv0\). Suppose instead that
\(\{v=0\}\ne\varnothing\) and \(L:=\sup_{\{v>0\}}|\nabla v|>1\). Following the
nearest-zero rescaling of \cite[proof of Proposition~A.5]{KW}, choose
\(x_k\in\{v>0\}\) with \(|\nabla v(x_k)|\to L\) and set
\(d_k=\operatorname{dist}(x_k,\{v=0\})>0\). The rescaled functions
\(w_k(x)=d_k^{-1}v(x_k+d_k x)\) are uniformly \(L\)-Lipschitz, vanish at a
point of \(\partial B_1\), and are harmonic on \(B_1\). Along a subsequence,
\(w_k\to w\) locally uniformly on \(\R^n\) and smoothly on \(B_1\), with
\(|\nabla w(0)|=L>0\) (so \(w>0\) on \(B_1\) by the minimum principle) and
\(|\nabla w|\le L\) on \(\{w>0\}\). Since \(\Delta|\nabla w|^2=2|D^2w|^2\ge0\)
on the connected component of \(\{w>0\}\) containing \(0\), the strong maximum
principle forces \(D^2w\equiv0\) there, making that component an open
half-space on which \(w\) has slope \(L\). On the other hand, selecting a
diagonal sequence \(v_{j(k)}\) approximating \(v\) within \(d_k/k\) on
\(x_k+d_k B_k\) exhibits \(w\) as a locally uniform limit of rescaled classical
solutions, so \(w\) is a viscosity supersolution in the sense of
\cite[Definition~A.1 and Lemma~A.3]{KW}; testing its half-space boundary with
an interior tangent ball forces \(L\le1\), a contradiction.$\square$
\end{proof}

\vskip0.23in

We next show that under a positive-phase Hessian bound, approximating free
boundaries decompose locally into either a single smooth graph or two ordered
oppositely oriented graphs enclosing a collapsing zero-phase layer (compare
\cite[proof of Lemma~5.1]{CFFS}).

\vskip0.23in

\begin{lemma}\label{lem:local-sheets}
Suppose classical solutions \(v_j\) converge locally uniformly near a point
\(q\) with uniformly bounded positive-phase Hessians there, and assume \(q\) is
an accumulation point of \(\FB(v_j)\). After passing to a subsequence, in a
fixed coordinate cylinder around \(q\) the free boundaries \(\FB(v_j)\) consist
either of a single smooth graph or of two ordered smooth graphs with opposite
inward normals separated by the zero phase. In smaller graph coordinates, the
graphs and their positive-side solutions converge smoothly; in the single-graph
case \(v\) has a nonempty open zero phase in the cylinder. In the
two-graph case the two limiting graphs meet at \(q\).
\end{lemma}

\begin{proof}
For \(q_j\in\FB(v_j)\) near \(q\), choose coordinates with vertical axis along
the inward unit normal \(e_j=\nabla v_j(q_j)\). Differentiating \(v_j=0\) along
\(\FB(v_j)\) with \(|\nabla v_j|=1\) bounds the second fundamental form of
\(\FB(v_j)\) by the positive-phase Hessian bound \(K\). Lifting radial rays
from the tangent disk of radius \(r\) via orthogonal projection keeps the unit
normal within angle \(\pi/4\) of \(e_j\) up to a uniform radius \(r\ll K^{-1}\),
yielding a smooth free-boundary graph over the \(r\)-disk with slope at most
\(CKr\). For a point \(a\) on the graph, \(\partial_{e_j}v_j(a)\ge 3/4\).
Integrating in the fixed vertical direction, up to any first exit, gives
\(v_j(a+te_j)\ge 3t/4-Kt^2/2>0\) for \(0<t\le cr\) and \(r\) small.
This rules out such an exit and gives a positive strip of height
comparable to \(r\) over a smaller disk. Retaining \(v_j\) above the
graph and setting it to zero below defines a classical one-sided solution in a
ball \(B_{cr}(q_j)\) that deviates from the unit half-space profile by at most
\(CKr^2\). Shrinking \(r\) allows us to apply one-sided flatness regularity
\cite{DeSilvaFlatness} in the quantitative form of \cite[Lemma~3.7]{CFFS}.
Together with implicit differentiation, this gives uniform
\(C^k\) bounds for every \(k\ge1\) on both the graph and its positive-side
solution.

\vskip0.23in

To count sheets near \(q\), extract a first limiting graph through a sequence
\(q_j\to q\) with limiting inward normal \(e\). Any second sequence of
free-boundary points tending to \(q\) off this sheet generates a second graph
whose limiting inward normal \(e'\) cannot be transverse to \(e\) (otherwise
its tangent plane at \(q\) would enter the open positive strip of the first
sheet, forcing free-boundary points of \(v_j\) inside \(\{v_j>0\}\)) and cannot
equal \(+e\) (otherwise one graph would lie strictly above the other over a
common disk and enter its positive strip, since \(\FB(v_j)\) is embedded). Thus
at most one additional sheet can approach \(q\), necessarily with opposite
limiting inward normal \(e'=-e\). If no additional sheet approaches \(q\),
a smaller cylinder contains only the first sheet for large \(j\). Otherwise
choose the second sheet through a sequence tending to \(q\), then shrink the
cylinder to exclude any further components by the same argument. In the
two-sheet case, non-intersection with each other's positive strips forces the
upper graph to have upward positive phase and the lower graph to have downward
positive phase. The region between the graphs lies in \(\{v_j=0\}\), and the
two limit graphs meet at \(q\) by construction.$\square$
\end{proof}

\vskip0.23in

\begin{lemma}\label{lem:sheet-contact}
Let \(v\in C(\R^n)\) be nonnegative, harmonic on \(\{v>0\}\), and satisfy
\(|\nabla v|\le1\) on \(\{v>0\}\). If two smooth one-sided Bernoulli free-boundary
graphs of \(v\) with opposite inward normals touch tangentially at a point and
bound opposite positive components, then \(v\) is globally a unit wedge.
\end{lemma}

\begin{proof}
On each adjacent positive branch, the function \(h=1-|\nabla v|^2\ge0\) is
superharmonic (\(\Delta h=-2|D^2v|^2\le0\)) and vanishes along its smooth free
boundary, where harmonicity and \(|\nabla v|=1\) give
\(\partial_{\nu_{\mathrm{in}}}h=2H\) for the inward mean curvature
\(H=\diver_{\FB}\nu_{\mathrm{in}}\). By the Hopf boundary point lemma, both
inward mean curvatures \(H_\pm\) are nonnegative. Writing the ordered graphs in
common tangent coordinates at the contact point as \(x_n=f_+(x')\ge f_-(x')\),
optimality at contact gives \(D^2(f_+-f_-)\ge0\) there, so
\(H_++H_-=-\Delta_{x'}(f_+-f_-)\le0\). Hence \(H_+=H_-=0\) at contact; the
strict Hopf lemma forces \(h\equiv0\), hence \(D^2v\equiv0\), on each local
positive branch. Unique continuation extends affineness across both positive
components in \(\R^n\), which by continuity and vanishing on their common
boundary hyperplane must be complementary unit half-spaces.$\square$
\end{proof}

\vskip0.23in

\begin{lemma}\label{lem:compact-hessian}
Let \(v_j\) be classical solutions on expanding balls \(D_j\uparrow\R^n\) with
\(0\in\FB(v_j)\), a uniform Lipschitz bound \(\|\nabla v_j\|_{L^\infty(D_j)}\le C_0\),
and uniformly bounded positive-phase Hessians on every compact subset of \(\R^n\). Assume there exists a fixed ball
on which \(\sup |D^2v_j|\ge c_0>0\) (including positive-side boundary traces).
Then, up to a subsequence, \(v_j\) converges locally uniformly on \(\R^n\) to a
nonflat entire classical solution \(v\), with smooth convergence in local
positive-side graph charts, \(\one_{\{v_j>0\}}\to\one_{\{v>0\}}\) in
\(L^1_{\mathrm{loc}}(\R^n)\), and
\[
  \operatorname{ind}(v)\le\liminf_{j\to\infty}\operatorname{ind}(v_j;D_j).
\]
\end{lemma}

\begin{proof}
Pass first to a subsequence realizing the index liminf. The
Arzel\'a--Ascoli theorem gives a locally uniform limit \(v\) after further
extraction. It is harmonic on \(\{v>0\}\) and satisfies \(|\nabla v|\le1\)
there by Lemma~\ref{lem:entire-slope}.

\vskip0.23in

We first rule out vanishing of the Hessian throughout \(\{v>0\}\).
If it vanished there, Lemma~\ref{lem:hessian-convergence} would give
\(\int_{B_{2R}\cap\{v_j>0\}}|D^2v_j|^n\dd x\to0\) for any fixed \(R\).
Choose \(B_R\) to contain the fixed ball from the peak hypothesis.
Hessian \(\varepsilon\)-regularity \cite[Lemma~3.4]{FRS}, rescaled from
\(B_{2R}\), would then give
\(\|D^2v_j\|_{L^\infty(B_R\cap\{v_j>0\})}\to0\). The same bound holds
for the one-sided traces, contradicting the peak hypothesis.

\vskip0.23in

Every \(q\in\FB(v)\) is approached by points of \(\FB(v_j)\). Otherwise
a fixed ball about \(q\) would contain no approximating free boundary
along a subsequence. In that ball each \(v_j\) is harmonic, including
the identically zero case, so the nonnegative harmonic limit vanishes
identically there by the minimum principle. Conversely, an accumulation
point of \(\FB(v_j)\) lies in \(\{v=0\}\) and is approached by positive
points of \(v\), by passing to the limit in the inward-segment bound
\eqref{eq:inward-segment} after extraction of the normals. Hence it
lies in \(\FB(v)\).

\vskip0.23in

Apply Lemma~\ref{lem:local-sheets} at a point of \(\FB(v)\). Two
limiting sheets would meet there and Lemma~\ref{lem:sheet-contact} would
force a unit wedge, excluded by nonvanishing of the positive-phase
Hessian. Thus the limit has one classical sheet and an open zero side
at each free boundary point. The local graph convergence follows by
the same compactness argument on every subsequence. In particular,
\(v\) is an entire classical solution whose Hessian is nonzero somewhere in \(\{v>0\}\).

\vskip0.23in

Positivity persists on compact subsets of \(\{v>0\}\). Each component
of \(\operatorname{int}\{v=0\}\) meets the zero side of a one-sheeted
boundary chart, since it has a boundary point in \(\FB(v)\).
A point of the component can be joined to such a zero side by a
compact path in the component. A small neighborhood of this path
avoids \(\FB(v_j)\) for large \(j\), by the characterization of its
accumulation points, and contains a point below the approximating
sheet in the zero phase. The zero phase therefore persists near
the original point. Since \(\mathcal H^n(\FB(v))=0\), this proves
\(\one_{\{v_j>0\}}\to\one_{\{v>0\}}\) in \(L^1_{\mathrm{loc}}(\R^n)\).

\vskip0.23in

Finally, given an \(m\)-dimensional subspace of \(C_c^\infty(\R^n)\) on which
\(\mathcal Q_v\) is negative definite, restricting a basis
\(\phi_1,\ldots,\phi_m\) to \(\overline{\{v_j>0\}}\) yields admissible tests on
\(D_j\) for large \(j\). Phase convergence in \(L^1_{\mathrm{loc}}\) and smooth
single-sheeted free-boundary convergence imply entrywise convergence of the
\(m\times m\) Gram matrices of \(\mathcal Q_{v_j}[\cdot;D_j]\) to the negative
definite matrix of \(\mathcal Q_v[\cdot;\R^n]\), proving lower semicontinuity
of the Morse index.$\square$
\end{proof}

\vskip0.23in

\subsection{Point-picking selection, interior curvature bounds, and finite-index entire profiles}\label{subsec:point-picking}

\begin{lemma}\label{lem:weighted-selection}
Let \(u_j\) be classical solutions on a bounded open set \(D\subset\R^n\) with
uniformly bounded Lipschitz norms on compact subsets of \(D\), and let \(d_j\)
be nonnegative \(1\)-Lipschitz functions on \(D\) satisfying
\(d_j\le\operatorname{dist}(\cdot,\partial D)\). If there exist points
\(x_j\in K\cap\{u_j>0\}\) in a fixed compact set \(K\Subset D\) such that
\(d_j(x_j)|D^2u_j(x_j)|\to\infty\), then there exist points \(y_j\) in a
compact subset of \(D\) and free-boundary centers \(z_j\in\FB(u_j)\) such that,
setting \(\lambda_j=|D^2u_j(y_j)|^{-1}\),
\[
  \lambda_j\to0,\qquad
  |y_j-z_j|=O(\lambda_j),\qquad
  \lambda_j=o\bigl(d_j(y_j)\bigr),
\]
and the rescalings \(\tilde u_j(x)=\lambda_j^{-1}u_j(z_j+\lambda_j x)\) have
positive-phase Hessians bounded by \(2\) on every fixed ball for large \(j\),
converging along a subsequence to a nonflat entire classical limit as in
Lemma~\ref{lem:compact-hessian}.
\end{lemma}

\begin{proof}
Set \(\rho_j=\frac14 d_j(x_j)\le\frac14\operatorname{dist}(x_j,\partial D)\).
For every \(y\in\overline B_{\rho_j}(x_j)\),
\[
  \operatorname{dist}(y,\partial D)
  \ge \tfrac34\operatorname{dist}(x_j,\partial D)
  \ge \tfrac34\operatorname{dist}(K,\partial D).
\]
Choose \(y_j\in\overline{\{u_j>0\}}\cap\overline{B_{\rho_j}(x_j)}\) maximizing
the weighted quantity \((\rho_j-|x-x_j|)|D^2u_j(x)|\) (including positive-side
boundary traces). Put \(t_j=\rho_j-|y_j-x_j|>0\) and \(a_j=|D^2u_j(y_j)|=\lambda_j^{-1}\);
comparing at \(x_j\) gives \(a_j t_j\ge\rho_j|D^2u_j(x_j)|\to\infty\). Because
all selection balls stay uniformly away from \(\partial D\), a single Lipschitz
constant \(L\) bounds \(|\nabla u_j|\) on them. If \(y_j\in\{u_j>0\}\),
interior harmonic derivative estimates on the ball around \(y_j\) of radius
\(\frac12\min\{t_j,\operatorname{dist}(y_j,\FB(u_j))\}\) give
\[
  a_j\le \frac{CL}{\min\{t_j,\operatorname{dist}(y_j,\FB(u_j))\}},
\]
which in view of \(a_j t_j\to\infty\) yields a nearest point
\(z_j\in\FB(u_j)\) at distance \(O(\lambda_j)\) from \(y_j\) (taking
\(z_j=y_j\) if \(y_j\in\FB(u_j)\)). By choice of \(y_j\), we have
\(|D^2u_j|\le 2a_j\) throughout \(B_{t_j/2}(y_j)\cap\{u_j>0\}\), while
\(1\)-Lipschitz continuity of \(d_j\) gives
\(d_j(y_j)\ge d_j(x_j)-\rho_j=\frac34 d_j(x_j)\), whence
\(\lambda_j/d_j(y_j)\le\frac{1}{3a_j\rho_j}\to0\) and \(t_j/\lambda_j\to\infty\).
Since \(|y_j-z_j|=O(\lambda_j)\), for large \(j\) we have
\(B_{t_j/4}(z_j)\subset B_{t_j/2}(y_j)\). Thus the rescalings have
uniform Lipschitz and Hessian bounds on \(B_{t_j/(4\lambda_j)}\), whose
radii tend to infinity. The rescaled Hessian peak stays at bounded
distance from the origin, so Lemma~\ref{lem:compact-hessian} applies.
$\square$
\end{proof}

\vskip0.23in

\begin{proof}[Proof of Corollary~\ref{cor:curv}]
By translation and Bernoulli scaling we may assume \(z=0\in\FB(u)\) and
\(R=1\). If \(\sup_{B_{1/2}\cap\{u>0\}}|D^2u_j|\to\infty\) along a sequence of
classical stable solutions \(u_j\) in \(B_1\) with \(0\in\FB(u_j)\), then
Lemma~\ref{lem:local-lipschitz} supplies uniform interior Lipschitz bounds, and
applying Lemma~\ref{lem:weighted-selection} with
\(d_j(x)=\operatorname{dist}(x,\partial B_1)\) on \(K=\overline{B_{1/2}}\)
produces a nonflat entire classical limit \(v\) with
\(\operatorname{ind}(v)\le 0\), contradicting Corollary~\ref{cor:classification}.
Finally, differentiating \(u=0\) twice tangentially along \(\FB(u)\) and using
\(|\nabla u|=1\) identifies \(\mathrm{II}\) with the tangential restriction of
the positive-side trace of \(D^2u\).$\square$
\end{proof}

\vskip0.23in

Combining Corollary~\ref{cor:classification} with \cite[Theorem~1.2]{FFSS} now
yields the complete list of nonflat finite-Morse-index bubbles in dimensions
three through six.

\vskip0.23in

\begin{corollary}\label{cor:finite}
For \(3\le n\le6\), every entire classical solution \(u\) to
\eqref{eq:bernoulli} of finite Morse index is axially symmetric up to a rigid
motion. If \(u\) is not one-dimensional, then up to rigid motion and Bernoulli
rescaling it is either the radial exterior-ball solution
\[
  u(x)=\frac{1-|x|^{2-n}}{n-2}\qquad (|x|>1),
\]
extended by zero on \(\{|x|\le1\}\), or an axially symmetric \emph{catenoid-type}
(two-ended) solution as in \cite[Theorem~1.2]{FFSS}, which is even in \(x_n\)
and satisfies \(\partial_{x_n}u>0\) on \(\{x_n>0\}\cap\{u>0\}\) and
\(\partial_r u<0\) on \(\{r>0\}\cap\{u>0\}\), where \(r=|x'|\).
Its zero set has the form
\[
  \bigl\{(x',x_n): r\ge1,\ |x_n|\le g(r)\bigr\},
  \qquad g(1)=0,\quad g'(r)>0\ (r>1),\quad \lim_{r\to\infty}g(r)<\infty.
\]
\end{corollary}

\begin{proof}
Theorem~1.2 of \cite{FFSS} applies for \(3\le n<\min\{n_{\mathrm{entire}}^*,n_{\mathrm{cone}}^*\}\)
(see \cite[(1.3)--(1.5)]{FFSS}). Since \(n_{\mathrm{entire}}^*=7\) by
Corollary~\ref{cor:critical} and \(n_{\mathrm{cone}}^*=7\) by the cone
classification \cite[Theorem~3 and Remark~4]{survey} and the seven-dimensional
example \cite{DSJ}, both thresholds equal \(7\), so the
classification holds for all \(3\le n\le6\).$\square$
\end{proof}

\vskip0.33in

\section{\bf Hessian energy identity at bounded Morse index}\label{sec:index-concentration}

To prove Theorem~\ref{thm:local-bubble-tree}, extracting scale-separated
bubbles via Lemma~\ref{lem:weighted-selection} is only the first step: one must
still show that no diffuse \(L^n\)-Hessian energy remains in the annular neck
regions separating bubble scales from macroscopic scales or from colliding
companion bubbles. We accomplish this in three steps: first pinning the Weiss
density across intermediate scales to force full asymptotic positive-phase
measure and classify intermediate limits; next establishing integrable Hessian
decay on intermediate annuli; and finally organizing multi-bubble collisions
via a Voronoi nearest-center scale tree.

\vskip0.23in

\subsection{Logarithmic Weiss bound and classification of intermediate blow-ups}\label{subsec:weiss}

Throughout this section we write \(A(a,b)=B_b\setminus\overline{B_a}\) for open
concentric annuli and \(\omega_n=|B_1|\) for the volume of the Euclidean unit
ball. For a function \(u\) on \(B_R\) with \(0\in\FB(u)\), recall the
scale-invariant Weiss energy functional
\[
  W(u,r)=r^{-n}\int_{B_r}\bigl(|\nabla u|^2+\one_{\{u>0\}}\bigr)\dd x
         -r^{-n-1}\int_{\partial B_r}u^2\dd\mathcal H^{n-1}.
\]
Our starting point is a quantitative logarithmic upper bound showing that any
Lipschitz classical solution cannot exceed the full-density threshold
\(\omega_n\) by more than \(O(1/\log(R/r))\).

\begin{lemma}\label{lem:weiss-cap}
Let \(u\) be a classical solution to \eqref{eq:bernoulli} in \(B_R\subset\R^n\)
satisfying \(u(0)=0\) and \(\|\nabla u\|_{L^\infty(B_R)}\le L\). Then
\[
  W(u,r)\le \omega_n+\frac{n\omega_n L^2}{2\log(R/r)}
  \qquad\forall\,0<r<R.
\]
In particular, every entire classical solution with \(0\in\FB(u)\) satisfies
\(W(u,r)\le\omega_n\) for all \(r>0\).
\end{lemma}

\begin{proof}
Integrating \(\Delta u=0\) against \(u\) by parts on \(B_r\cap\{u>0\}\) and
using Weiss's monotonicity identity for stationary solutions
\cite[Theorem~3.1]{WeissWeak} (see also
\cite[Proposition~3.4]{KriventsovWeiss}) give, for
almost every \(r\in(0,R)\),
\begin{align*}
  W(u,r)&=r^{-n}\bigl|B_r\cap\{u>0\}\bigr|
          +r^{-n}\int_{\partial B_r}u\bigl(\partial_r u-u/r\bigr)\dd\mathcal H^{n-1},\\
  \partial_r W(u,r)&=2r^{-n}\int_{\partial B_r}\bigl(\partial_r u-u/r\bigr)^2\dd\mathcal H^{n-1}.
\end{align*}
Bounding the first term in \(W(u,r)\) by \(\omega_n\), using \(0\le u\le Lr\)
on \(\partial B_r\) with \(\mathcal H^{n-1}(\partial B_r)=n\omega_n r^{n-1}\),
and applying the Cauchy--Schwarz inequality yield the Riccati-type differential
inequality
\[
  \bigl(W(u,r)-\omega_n\bigr)_+^2
  \le \tfrac12 n\omega_n L^2 r\,\partial_r W(u,r).
\]
Whenever \(W(u,r)>\omega_n\), monotonicity guarantees \(W(u,s)\ge W(u,r)>\omega_n\)
for all \(s\in[r,R)\); integrating \(-\partial_s(W(u,s)-\omega_n)^{-1}\) over
\([r,R']\) and letting \(R'\uparrow R\) gives the stated bound. For an entire
solution, \(|\nabla u|\le1\) globally by \cite[Lemma~3.5]{FRS}, so letting
\(R\to\infty\) gives \(W(u,r)\le\omega_n\).$\square$
\end{proof}

\vskip0.23in

\begin{proposition}\label{prop:intermediate-host}
Let \(3\le n\le6\), and let \(u_j\) be classical solutions in a fixed
neighborhood of an interior point \(p\) with uniformly bounded Lipschitz norms
and Morse indices at most \(I\). Suppose \(z_j\in\FB(u_j)\) converge to \(p\)
and, for scales \(\lambda_j\to0\), the bubble rescalings
\(\lambda_j^{-1}u_j(z_j+\lambda_j x)\) converge as in
Lemma~\ref{lem:compact-hessian} to a nonflat entire finite-index solution
\(V\). Then for every intermediate scale sequence \(s_j\to0\) satisfying
\(\lambda_j/s_j\to0\), every locally uniform subsequential limit \(w\) of
\(w_j(x)=s_j^{-1}u_j(z_j+s_j x)\) is either identically zero or a centered unit
wedge \(|x\cdot e|\), and along the same subsequence
\(\one_{\{w_j>0\}}\to 1\) strongly in \(L^1_{\mathrm{loc}}(\R^n)\).
\end{proposition}

\begin{proof}
Pass to a subsequence along which \(w_j\to w\) locally uniformly on \(\R^n\).
We first note that for any locally uniform limit \(f\) of uniformly Lipschitz
classical solutions, \(f\) is nonnegative, locally Lipschitz, and harmonic on
\(\{f>0\}\); testing \(\Delta f=0\) against \(\psi(f-\varepsilon)_+\) for
\(\psi\in C_c^\infty(\R^n)\) and letting \(\varepsilon\downarrow0\) (using
\(\nabla f=0\) a.e.\ on \(\{f=0\}\)) gives the distributional identity
\[
  \int_{\R^n}\psi|\nabla f|^2\dd x
  =-\int_{\R^n}f\nabla f\cdot\nabla\psi\dd x
  =\tfrac12\int_{\R^n}f^2\Delta\psi\dd x.
\]
Because the same identity holds for each \(w_j\), uniform convergence implies
convergence of localized \(L^2\) gradient norms, upgrading weak \(H^1_{\mathrm{loc}}\)
convergence to strong convergence \(w_j\to w\) in \(H^1_{\mathrm{loc}}(\R^n)\).

\vskip0.23in

Next we compute \(\lim_{R\to\infty}W(V,R)\). By Corollary~\ref{cor:finite}, the
blow-down \(V_R(x)=R^{-1}V(Rx)\) as \(R\to\infty\) converges locally uniformly
(and hence strongly in \(H^1_{\mathrm{loc}}\)) to \(0\) when \(V\) is an
exterior-ball solution, and to a centered unit wedge \(|x\cdot e|\) when \(V\)
is a catenoid-type solution (via the asymptotic end expansion of
\cite[Remark~1.4]{FFSS}). In both cases \(\one_{\{V_R>0\}}\to 1\) in
\(L^1_{\mathrm{loc}}(\R^n)\), and the Dirichlet and boundary terms in
\(W(V,R)=W(V_R,1)\) cancel in the limit (trivially for \(0\), and by
integration by parts on the two complementary half-spaces for \(|x\cdot e|\)),
giving \(\lim_{R\to\infty}W(V,R)=\omega_n\).

Now fix any \(t>0\) and \(R>0\). Since \(s_j t>\lambda_j R\) for large \(j\),
Weiss monotonicity and smooth/phase convergence at the bubble scale give
\[
  \liminf_{j\to\infty}W(w_j,t)
  =\liminf_{j\to\infty}W(u_j(z_j+\cdot),s_j t)
  \ge \lim_{j\to\infty}W(u_j(z_j+\cdot),\lambda_j R)
  =W(V,R).
\]
Letting \(R\to\infty\) yields \(\liminf_{j\to\infty}W(w_j,t)\ge\omega_n\),
while applying Lemma~\ref{lem:weiss-cap} from radius \(s_j t\) up to a fixed
macroscopic ball around \(p\) gives \(\limsup_{j\to\infty}W(w_j,t)\le\omega_n\).
Hence \(W(w_j,t)\to\omega_n\) for every \(t>0\). Passing to the limit via
strong \(H^1_{\mathrm{loc}}\) convergence in the integrated Weiss identity
\[
  W(w_j,b)-W(w_j,a)
  =2\int_{A(a,b)}|x|^{-n}\bigl(\partial_r w_j-w_j/|x|\bigr)^2\dd x
  \qquad (0<a<b<\infty)
\]
shows that \(w\) is homogeneous of degree one. Furthermore, if \(\chi\in L^\infty(\R^n)\)
is any weak-\(*\) subsequential limit of \(\one_{\{w_j>0\}}\), then passing to
the limit in \(W(w_j,t)\) gives
\[
  \omega_n
  =t^{-n}\int_{B_t}\bigl(|\nabla w|^2+\chi\bigr)\dd x
   -t^{-n-1}\int_{\partial B_t}w^2\dd\mathcal H^{n-1}.
\]
Because \(w\) is \(1\)-homogeneous and \(\Delta(w^2)=2|\nabla w|^2\)
distributionally, \(\int_{B_t}|\nabla w|^2\dd x=t^{-1}\int_{\partial B_t}w^2\dd\mathcal H^{n-1}\)
for almost every \(t>0\), leaving \(\int_{B_t}\chi\dd x=\omega_n t^n=|B_t|\).
Since \(0\le\chi\le1\), we conclude \(\chi\equiv1\) and therefore
\(\one_{\{w_j>0\}}\to1\) strongly in \(L^1_{\mathrm{loc}}(\R^n)\).

\vskip0.23in

It remains to classify the \(1\)-homogeneous limit \(w\). Diagonalizing over
all balls with rational centers and rational radii, pass to a further
subsequence along which each rational ball is either eventually stable or
eventually unstable for \(w_j\). The complement \(\Sigma\) of the union
\(\mathcal O\) of eventually stable rational balls contains at most \(I\)
points (as \(I+1\) distinct points in \(\Sigma\) could be separated by \(I+1\)
disjoint eventually unstable rational balls, producing \(I+1\) mutually
orthogonal negative-definite directions). On every compact subset of
\(\mathcal O\), Corollary~\ref{cor:curv} and the interior harmonic
estimate give uniform positive-phase Hessian bounds for \(w_j\).
Indeed, if \(B_{3r}(x)\) is contained in an eventually stable rational
ball, use the curvature estimate centered at a point of
\(\FB(w_j)\cap B_{r/2}(x)\), when such a point exists, with radius
\(2r\). Otherwise \(w_j\) is harmonic in \(B_{r/2}(x)\), including the
identically zero case. This bounds the Hessian in \(B_{r/4}(x)\) in
either case. Lemma~\ref{lem:entire-slope} gives \(|\nabla w|\le1\)
on \(\{w>0\}\).

\vskip0.23in

Suppose \(w\not\equiv0\). Then \(\FB(w)\setminus\{0\}\) must be nonempty:
otherwise \(w>0\) would be harmonic on \(\R^n\setminus\{0\}\), which for
\(n\ge3\) has a removable singularity at \(0\) and would force \(w\) to be a
positive constant, incompatible with \(1\)-homogeneity. By \(1\)-homogeneity,
\(\FB(w)\setminus\{0\}\) contains an entire open ray, which necessarily meets
\(\mathcal O=\R^n\setminus\Sigma\) since \(\#\Sigma\le I\). Choose
\(q\in\FB(w)\cap\mathcal O\) on that ray. If a ball about \(q\) contained
no approximating free boundary along a subsequence, its nonnegative harmonic
limit would vanish identically there by the minimum principle, contrary to
\(q\in\FB(w)\). Thus Lemma~\ref{lem:local-sheets} applies at \(q\).
A single limiting sheet would leave an open zero side. Two sheets separated
at \(q\) would leave an open zero layer. In either case a smaller ball would
remain in the zero phase of \(w_j\), contrary to
\(\one_{\{w_j>0\}}\to1\) in \(L^1_{\mathrm{loc}}\). The two limit sheets
therefore meet at \(q\), and Lemma~\ref{lem:sheet-contact} identifies \(w\)
with a unit wedge, centered at the origin by \(1\)-homogeneity.$\square$
\end{proof}

\vskip0.23in

\subsection{Annular branch splitting and vanishing of neck energy}\label{subsec:annuli}

When an intermediate profile is a centered unit wedge rather than zero, the
positivity set on a long dyadic annulus surrounds a thin collapsing zero-phase
sheet. We first show that local two-graph descriptions across consecutive
scales glue coherently into two globally disjoint classical branches trapped
between parallel shifted half-spaces.

\vskip0.23in

\begin{lemma}\label{lem:glue-neck-branches}
Let \(0<s_0<s_1/16\), and let \(u\) be a classical solution to
\eqref{eq:bernoulli} on \(A(s_0/4,16s_1)\). Assume that for every
\(s\in[s_0,s_1]\) there exists \(e_s\in\mathbb S^{n-1}\) satisfying
\[
  \bigl|u(x)-|e_s\cdot x|\bigr|\le\eta s
  \qquad\forall\,x\in A(s/4,16s),
\]
such that in orthogonal coordinates \(x=\xi+\tau e_s\) (\(\xi\perp e_s\)), for
every \(\xi\) with \(\frac38 s<|\xi|<9s\) the intersection of \(\{u=0\}\) with
\(\{\xi+\tau e_s:|\tau|<s/16\}\) is a closed interval whose two endpoints vary
smoothly with \(\xi\) and have absolute height at most \(\eta s\). Provided
\(\eta>0\) is smaller than a dimensional constant, on the trimmed annulus
\(A(\rho,R)\) with \(\rho=4s_0\) and \(R=s_1/4\) one can decompose
\(u=v_++v_-\) into two classical solutions with disjoint positivity sets, and
for every \(2\rho<s<R/2\) orient \(e_s\) so that
\[
  (\pm e_s\cdot x-C\eta s)_+\le v_\pm(x)\le (\pm e_s\cdot x+C\eta s)_+
  \qquad\forall\,x\in A(s/2,2s),
\]
with \(C>0\) independent of the scale ratio \(s_1/s_0\).
\end{lemma}

\begin{proof}
At each scale \(s\in[s_0,s_1]\), the two ordered endpoint graphs partition
\(A(\frac{7}{16}s,8s)\cap\{u>0\}\) into two locally constant sign classes
(determined by \(\operatorname{sgn}\tau\) outside \(|\tau|\le\eta s\), and
above or below the zero-phase interval inside \(|\tau|\le\eta s\), where
\(|\xi|>\frac38 s\) holds automatically). Consider the dyadic ladder
\(t=2^k s_0\le s_1\) with retained shells \(A(t,4t)\). Evaluating the wedge
approximations at scales \(t\) and \(2t\) along the equator
\(\{x\perp e_t:|x|=3t\}\) gives \(|P_{e_t^\perp}e_{2t}|\le C\eta\), allowing
successive orientation of the dyadic normals so that \(|e_t-e_{2t}|\le C\eta\).
On the overlap \(A(2t,4t)\), the two scale labels manifestly agree in the deep
caps \(\{\pm e_t\cdot x>C_0\eta t\}\); moreover, any positive point of
\(A(2t,4t)\) outside those caps can be translated by at most \(C_1\eta t\)
along \(\pm e_t\) away from the zero interval into its own deep cap along a
segment that remains strictly positive and inside both extended shells
\(A(\frac{7}{16}t,8t)\) and \(A(\frac{7}{8}t,16t)\), forcing the \(t\)- and
\(2t\)-labels to agree throughout \(A(2t,4t)\cap\{u>0\}\).

\vskip0.23in

Setting \(v_\pm=u\) on the resulting global positive components in \(A(\rho,R)\)
and zero elsewhere defines two disjoint classical solutions on \(A(\rho,R)\),
since locally each inherits one of the two smooth one-sided sheets of \(u\).
Finally, for any intermediate \(s\in(2\rho,R/2)\), choosing the dyadic scale in
\([s/2,s]\) orients \(e_s\) within \(O(\eta)\) of the overlapping dyadic axes,
identifies the local \(s\)-branches with \(v_\pm\) on \(A(s/2,2s)\), and yields
the two-sided half-space trapping bounds both in the caps
\(\{\pm e_s\cdot x>\eta s\}\) (where \(v_\pm=u\) and \(v_\mp=0\)) and in the
thin strip \(\{|e_s\cdot x|\le\eta s\}\) (where \(0\le v_\pm\le 2\eta s\)).$\square$
\end{proof}

\vskip0.23in

\begin{lemma}\label{lem:abstract-neck}
Let \(3\le n\le6\), let \(a_j/b_j\to0\), and suppose classical solutions \(u_j\)
centered at \(z_j\in\FB(u_j)\) have uniformly bounded Lipschitz norms on
\(B_{16b_j}(z_j)\). Assume that for every sequence \(s_j\) with
\(a_j/s_j\to0\) and \(s_j/b_j\to0\), the rescalings
\(s_j^{-1}u_j(z_j+s_j x)\) have uniformly bounded positive-phase Hessians on
compact annuli in \(\R^n\setminus\{0\}\) and admit only \(0\) or centered unit
wedges as locally uniform subsequential limits. Then
\begin{equation}\label{eq:abstract-neck-energy}
  \lim_{M\to\infty}\limsup_{j\to\infty}
  \mu_{u_j}\bigl(B_{b_j/M}(z_j)\setminus B_{Ma_j}(z_j)\bigr)=0.
\end{equation}
\end{lemma}

\begin{proof}
Translate \(z_j\) to the origin. Choose \(\eta>0\) small enough for
Lemma~\ref{lem:glue-neck-branches} and so that its trapping constant
\(C\eta\) meets the threshold in \cite[Theorem~2.2]{FFSS} for exponent
\(1/2\). We also choose \(\eta\) small enough for the one-sided
\(\varepsilon\)-regularity estimate below. We claim that there exists a fixed
\(S\ge2\) such that for all large \(j\) and every \(s\in[Sa_j,b_j/S]\), one of
two mutually exclusive regimes holds:
\begin{enumerate}
\item[\textup{(1)}] \(\sup_{A(s,2s)}u_j\le\eta s\) and \(\FB(u_j)\cap A(s,2s)=\varnothing\);
\item[\textup{(2)}] the two-graph wedge hypotheses of
Lemma~\ref{lem:glue-neck-branches} hold at scale \(s\).
\end{enumerate}
Were this false, there would exist \(S_k\to\infty\), \(j_k\to\infty\), and
scales \(s_k\in[S_k a_{j_k},b_{j_k}/S_k]\) violating both~(1) and~(2). Passing
to a locally uniform limit of \(w_k(x)=s_k^{-1}u_{j_k}(s_k x)\), if
\(w_k\to0\), then no free-boundary points can persist in \(A(1,2)\) by the
inward positivity bound \eqref{eq:inward-segment}, so~(1) holds for large \(k\).
If instead \(w_k\to|e\cdot x|\), then \(\FB(w_k)\) on compact annuli concentrates
inside a vanishing strip around \(e^\perp\). For any convergent
sequence of normals with limit \(\nu\), pass to the limit in the first
bound of \eqref{eq:inward-segment} at a fixed sufficiently small
\(t>0\) to get \(t|\nu\cdot e|\ge t-Kt^2/2\), with \(K\) the local
uniform Hessian bound. Letting \(t\downarrow0\) shows
\(|\nu_k\cdot e|\to1\) uniformly in the annular region in question.
Along each vertical segment
\(\{\xi+\tau e:|\tau|\le 1/16\}\) with \(\frac38<|\xi|<9\), interior harmonic
convergence gives \(\partial_e w_k(\xi\pm\frac{1}{16}e)\to\pm1\); the uniform
Hessian bound rules out both the absence of zeros (which would force
\(\partial_e w_k\) to jump by nearly \(2\) across \([-\delta,\delta]\) for
\(\delta\ll C^{-1}\)) and the presence of multiple zero intervals (which would
trap a vanishingly short positive interval with opposite unit endpoint slopes),
while transversality and one-sidedness produce a nondegenerate zero interval
whose endpoints form smooth graphs of height at most \(\eta\), fulfilling~(2)
for large \(k\).

\vskip0.23in

Because the continuous scale function \(s\mapsto s^{-1}\sup_{A(s,2s)}u_j\) lies
in \([0,\eta]\) under~(1) and in \([1-\eta,\infty)\) under~(2), connectedness
of \([Sa_j,b_j/S]\) forces a single regime to hold across the entire interval
for each large \(j\). Put \(\rho=4Sa_j\) and \(R=b_j/(4S)\).

\vskip0.23in

In regime~(1), \(u_j\) is harmonic (possibly identically zero) on \(A(\rho,R)\).
Decomposing \(g=\partial_m u_j\) on \(A(2\rho,R/2)\) via the harmonic extension
\(G\) of \(g|_{\partial B_{R/2}}\) to \(B_{R/2}\), the maximum principle gives
\(|g-G|\le 2L(2\rho/|x|)^{n-2}\), and interior harmonic derivative bounds on
balls of radius \(\asymp s\) around \(|x|=s\in(8\rho,R/8)\) yield
\[
  s|D^2u_j(x)|\le CL\bigl[(\rho/s)^{n-2}+s/R\bigr].
\]
In regime~(2), Lemma~\ref{lem:glue-neck-branches} splits \(u_j=v_++v_-\) on
\(A(\rho,R)\), and applying the annular improvement-of-flatness theorem
\cite[Theorem~2.2]{FFSS} with decay exponent \(1/2\) to each \(v_\pm\) supplies
unit-slope affine approximations on \(A(s/2,2s)\cap\{v_\pm>0\}\) with height
error
\[
  h\le C\eta s\bigl[(\rho/s)^{n-3/2}+(s/R)^{1/2}\bigr].
\]
At a point of \(\{|x|=s\}\cap\{v_\pm>0\}\) at distance at least
\(s/32\) from \(\FB(v_\pm)\), apply the interior harmonic estimate
to the difference from the affine fit \(L\). Otherwise take a nearest
free-boundary point \(q\). Continuity gives \(|L(q)|\le h\), so
\(|v_\pm-(L-L(q))|\le2h\) on the positive phase of
\(B_{s/8}(q)\subset A(s/2,2s)\). Translate to \(q\), rotate the unit
slope, and rescale by \(s/8\) to apply \cite[Lemma~3.7]{CFFS}.
The rescaled flatness is at most \(16h/s=O(\eta)\), within the
required threshold by the choice of \(\eta\). Both estimates
give \(s|D^2v_\pm(x)|\le Ch/s\). Since \(n\ge3\), in
both regimes we have
\(s|D^2u_j(x)|\le C(L+\eta)\bigl[(\rho/s)^{1/2}+(s/R)^{1/2}\bigr]\) throughout
\(A(8\rho,R/8)\cap\{u_j>0\}\). Raising to the \(n\)-th power and integrating
against \(s^{-1}\dd s\) over \(Ma_j<|x|<b_j/M\) (for \(M\gg S\)) bounds the
annular Hessian mass by \(C(L+\eta)^n(S/M)^{n/2}\), which vanishes as
\(M\to\infty\).$\square$
\end{proof}

\vskip0.23in

The same exterior affine-fit estimate immediately yields finiteness of the
total Hessian energy for every finite-Morse-index bubble.

\vskip0.23in

\begin{corollary}\label{cor:bubble-finite-energy}
Let \(3\le n\le6\), and let \(V\) be an entire classical solution to
\eqref{eq:bernoulli} of finite Morse index. Then \(\mu_V(\R^n)<\infty\).
\end{corollary}

\begin{proof}
If \(V\) is one-dimensional, \(\mu_V(\R^n)=0\); if \(V\) is a radial
exterior-ball solution, \(|D^2V(x)|=O(|x|^{-n})\) as \(|x|\to\infty\), which is
integrable in \(L^n(\R^n\setminus B_1)\). For an axially symmetric
catenoid-type solution, \cite[Proposition~3.1]{FFSS} splits \(V\)
outside a fixed ball into two classical branches with disjoint
supports and centered annular heights \(o(s)\) as \(s\to\infty\).
Choose \(\rho\) large enough that both branches meet the smallness
hypothesis of \cite[Theorem~2.2]{FFSS} on \(A(\rho,R)\) for every
sufficiently large \(R\). Applying the theorem with exponent \(1/2\),
using the same one-sided estimate as above and letting \(R\to\infty\)
gives \(|x|\,|D^2V(x)|\le C(\rho/|x|)^{n-3/2}\) on the positive phase
for \(|x|>8\rho\). Thus \(|D^2V|^n\le C_1|x|^{-n(n-1/2)}\) there,
which is integrable at infinity.$\square$
\end{proof}

\vskip0.23in

\subsection{Voronoi scale tree and proof of Theorem~\ref{thm:local-bubble-tree}}\label{sec:local-bubble-tree}

\begin{proof}[Proof of Theorem~\ref{thm:local-bubble-tree}]
We organize the proof into four stages.

\smallskip
\noindent\emph{Stage 1: Iterative curvature selection and index accounting.}
By Lemma~\ref{lem:local-lipschitz} and \(\FB(u_j)\cap B_1\ne\varnothing\), the
sequence \(\{u_j\}\) is locally equi-Lipschitz in \(B_2\) and converges locally
uniformly after passing to a subsequence to a locally Lipschitz limit \(U\)
that is harmonic on \(\{U>0\}\). Suppose inductively that \(m\ge0\) bubble
sequences \((z_j^i,\lambda_j^i,V^i)_{1\le i\le m}\) have already been extracted,
and define the \(1\)-Lipschitz distance weight
\[
  d_j(x)=\min\Bigl\{\operatorname{dist}(x,\partial B_2),\
                    \min_{1\le i\le m}|x-z_j^i|\Bigr\}
\]
(omitting the inner minimum when \(m=0\)). If \(d_j|D^2u_j|\) is unbounded on
some compact subset of \(B_2\), Lemma~\ref{lem:weighted-selection} extracts a
new center \(z_j^{m+1}\in\FB(u_j)\) converging to some \(p^{m+1}\in B_2\) and a
scale \(\lambda_j^{m+1}\to0\) along a subsequence such that the rescalings
converge to a nonflat entire classical solution \(V^{m+1}\). By
Corollary~\ref{cor:classification}, \(\operatorname{ind}(V^{m+1})\ge1\), while
Lemma~\ref{lem:compact-hessian} gives \(\operatorname{ind}(V^{m+1})\le I\).
Moreover, for each \(1\le i\le m\), we have
\(|z_j^{m+1}-z_j^i|/\lambda_j^{m+1}\to\infty\) because
\(\lambda_j^{m+1}=o(d_j(y_j))\) and \(|y_j-z_j^{m+1}|=O(\lambda_j^{m+1})\);
conversely, if \(|z_j^{m+1}-z_j^i|/\lambda_j^i\) remained bounded, then
\(|y_j-z_j^i|\le 2|z_j^{m+1}-z_j^i|=O(\lambda_j^i)\) for large \(j\), and the
bounded Hessian of the \(i\)-th bubble rescaling would force
\(d_j(y_j)/\lambda_j^{m+1}\le |y_j-z_j^i|\,|D^2u_j(y_j)|\le C\), a
contradiction.

\vskip0.23in

Thus all extracted bubbles satisfy the separation in part~(ii).
Choosing for each \(1\le i\le m+1\) an \(\operatorname{ind}(V^i)\)-dimensional
subspace of admissible test functions supported in a fixed ball \(B_{R_i}(0)\)
on which \(\mathcal Q_{V^i}\) is negative definite, mutual scale separation
ensures that for all sufficiently large \(j\) the rescaled support balls
\(B_{R_i\lambda_j^i}(z_j^i)\) are pairwise disjoint subsets of \(B_2\).
Consequently, the transferred negative-definite subspaces have mutually
disjoint supports in \(B_2\), forcing \(\sum_{i=1}^{m+1}\operatorname{ind}(V^i)\le I\). Since each bubble consumes at least one
unit of index, the selection terminates at some \(0\le N\le I\) with the
uniform residual curvature bound
\begin{equation}\label{eq:cluster-residual-bound}
  |D^2u_j(x)|\le \frac{C_K}{d_j(x)}
  \qquad\forall\,x\in K\cap\{u_j>0\},\ d_j(x)>0,
\end{equation}
on every compact subset \(K\Subset B_2\). Part~(i) now follows from
Corollaries~\ref{cor:finite} and~\ref{cor:bubble-finite-energy}.

\vskip0.23in

\smallskip
\noindent\emph{Stage 2: Voronoi cell partition and scale hierarchy at a concentration point.}
Fix a limiting concentration point \(p\in\{p^1,\ldots,p^N\}\), put
\(J_p=\{1\le i\le N:p^i=p\}\), and choose \(b>0\) small enough that
\(B_{32b}(p)\Subset B_2\) are pairwise disjoint across distinct concentration
points. Partition \(\R^n\) into the Voronoi nearest-center cells \(E_j^i\)
associated with \(\{z_j^i:i\in J_p\}\) (breaking ties on bisecting hyperplanes
arbitrarily), so that
\[
  x\in E_j^i\quad\Longrightarrow\quad
  |x-z_j^k|\ge |x-z_j^i|\qquad\forall\,k\in J_p.
\]
For each fixed \(i\in J_p\), consider the finite list of scales
\(\{\lambda_j^i\}\cup\{|z_j^i-z_j^k|:k\in J_p\setminus\{i\}\}\cup\{b\}\). Passing
to a subsequence so that all pairwise scale ratios converge in \([0,\infty]\),
merge scales with finite nonzero limiting ratios and choose one representative
from each equivalence class to form a strictly ordered hierarchy
\[
  d_{j,0}\ll d_{j,1}\ll\cdots\ll d_{j,m},
  \qquad d_{j,0}=\lambda_j^i,\quad d_{j,m}=b,
\]
where \(a_j\ll b_j\) signifies \(a_j/b_j\to0\) (notice \(d_{j,0}\ll d_{j,1}\)
by part~(ii)). Take one further subsequence, common to all finitely many
groups and scales, on which every bounded vector sequence
\((z_j^k-z_j^i)/d_{j,h}\) with \(i,k\in J_p\) and \(0<h<m\) converges.
The norms of the other such vector sequences tend to infinity by the
choice of the scale ratios.

\vskip0.23in

\smallskip
\noindent\emph{Stage 3: Vanishing of gap and transition-band energies.}
First consider any scale gap \(d_{j,h}\ll s_j\ll d_{j,h+1}\) (\(0\le h<m\)).
Rescaling \(u_j\) around \(z_j^i\) by \(s_j\), every other center \(z_j^k\)
either collapses onto the origin (if \(|z_j^k-z_j^i|\lesssim d_{j,h}\)) or
escapes to infinity (if \(|z_j^k-z_j^i|\gtrsim d_{j,h+1}\)). The residual bound
\eqref{eq:cluster-residual-bound} therefore furnishes uniform Hessian bounds on
compact annuli in \(\R^n\setminus\{0\}\), and Proposition~\ref{prop:intermediate-host}
restricts every subsequential limit to \(0\) or a centered unit wedge. Applying
Lemma~\ref{lem:abstract-neck} shows that the full gap annulus
\(A(Md_{j,h},d_{j,h+1}/M)+z_j^i\) carries asymptotically vanishing
\(\mu_{u_j}\)-mass as \(j\to\infty\) followed by \(M\to\infty\).

\vskip0.23in

Next consider an intermediate representative scale \(d_{j,h}\) with \(0<h<m\),
and restrict the transition band \(A(d_{j,h}/M,Md_{j,h})+z_j^i\) to the Voronoi
cell \(E_j^i\). After rescaling around \(z_j^i\) by \(d_{j,h}\), every point
\(y\) in the rescaled restricted band satisfies \(|y|\ge 1/M\), and by the
Voronoi property its distance to every rescaled companion center
\((z_j^k-z_j^i)/d_{j,h}\) (\(k\in J_p\)) is also at least \(|y|\ge 1/M\), while
centers outside \(J_p\) escape to infinity. Thus, for each fixed
\(M\ge2\) and all large \(j\), the rescaled restricted band lies in
\(\{1/M\le |y|\le M\}\) at distance at least \(1/(2M)\) from each of
the finitely many limiting centers selected above. The resulting fixed
compact set has a neighborhood where
\eqref{eq:cluster-residual-bound} bounds the rescaled Hessians uniformly.
By Proposition~\ref{prop:intermediate-host}, every locally uniform
limit is \(0\) or a unit wedge, both with vanishing positive-phase
Hessian. Lemma~\ref{lem:hessian-convergence} with \(p=n\) therefore
gives, at each fixed \(M\),
\[
  \lim_{j\to\infty}
  \mu_{u_j}\bigl(E_j^i\cap(A(d_{j,h}/M,Md_{j,h})+z_j^i)\bigr)=0.
\] Meanwhile, on the core ball \(B_{M\lambda_j^i}(z_j^i)\)
---which lies entirely inside \(E_j^i\) for large \(j\) since
\(\lambda_j^i\ll |z_j^i-z_j^k|\)---Lemma~\ref{lem:hessian-convergence} gives
\(\lim_{j\to\infty}\mu_{u_j}(B_{M\lambda_j^i}(z_j^i))=\mu_{V^i}(B_M)\).

\smallskip
\noindent\emph{Stage 4: Measure synthesis.}
For fixed large \(M\) and all sufficiently large \(j\), decomposing each cell
\(E_j^i\cap B_{b/M}(z_j^i)\) radially into its core ball, restricted transition
bands, and gap annuli covers the union
\[
  \mathcal U_{j,M}=\bigcup_{i\in J_p}\bigl(E_j^i\cap B_{b/M}(z_j^i)\bigr),
\]
yielding two-sided mass bounds
\[
  \sum_{i\in J_p}\mu_{u_j}\bigl(B_{M\lambda_j^i}(z_j^i)\bigr)
  \le \mu_{u_j}(\mathcal U_{j,M})
  \le \sum_{i\in J_p}\mu_{u_j}\bigl(B_{M\lambda_j^i}(z_j^i)\bigr)+\mathcal R_{j,M},
  \qquad \lim_{M\to\infty}\limsup_{j\to\infty}\mathcal R_{j,M}=0.
\]
Since \(z_j^i\to p\) for all \(i\in J_p\), we have the geometric sandwich
\(B_{b/(2M)}(p)\subset\mathcal U_{j,M}\subset B_{2b/M}(p)\) for all large \(j\).
Evaluating the upper bound at a single fixed \(M\) gives uniform local mass
bounds for \(\{\mu_{u_j}\}\) near each \(p\), and away from the finite set
\(\{p^1,\ldots,p^N\}\), \eqref{eq:cluster-residual-bound} and
Lemma~\ref{lem:hessian-convergence} give strong \(L^1_{\mathrm{loc}}\)
convergence of \(|D^2u_j|^n\one_{\{u_j>0\}}\) to \(|D^2U|^n\one_{\{U>0\}}\)
(whose integrability across \(\{p^1,\ldots,p^N\}\) follows from Fatou's lemma
and the uniform local mass bound). Finally, sending \(j\to\infty\) and then
\(M\to\infty\) along continuity radii \(b/(2M)\) and \(2b/M\) of a weak-\(*\)
Radon limit identifies the atomic weight at \(p\) as
\(\sum_{i\in J_p}\lim_{M\to\infty}\mu_{V^i}(B_M)=\sum_{i\in J_p}\mu_{V^i}(\R^n)\),
establishing \eqref{eq:local-multi-energy}.$\square$
\end{proof}

\vskip0.23in

\begin{corollary}\label{cor:local-uniform-energy}
Under the hypotheses of Theorem~\ref{thm:local-bubble-tree}, for every
\(0<r<2\) there exists \(C=C(n,I,r)>0\) such that \(\mu_{u_j}(B_r)\le C\).
\end{corollary}

\begin{proof}
Were \(\mu_{u_j}(B_r)\to\infty\) along a sequence satisfying the hypotheses,
testing \eqref{eq:local-multi-energy} along a convergent subsequence against a
nonnegative function \(\chi\in C_c(B_2)\) with \(\chi\equiv1\) on
\(\overline{B_r}\) would yield a finite limit, a contradiction.$\square$
\end{proof}

\vskip0.23in

\begin{corollary}\label{cor:local-index-saturated}
Under the hypotheses of Theorem~\ref{thm:local-bubble-tree} with \(I\ge1\),
suppose there exist \(I\) sequences of positive-phase points
\(x_j^\ell\to p^\ell\) (\(\ell=1,\ldots,I\)) with pairwise distinct limits
\(p^1,\ldots,p^I\in B_2\) along which \(|D^2u_j(x_j^\ell)|\to\infty\). Then on
the subsequence in Theorem~\ref{thm:local-bubble-tree} precisely one bubble forms at each \(p^\ell\),
each bubble has Morse index one, no other atoms occur, and for \(I=1\) the
unique atom carries the full Hessian energy of its bubble.
\end{corollary}

\begin{proof}
By \eqref{eq:cluster-residual-bound}, curvature blow-up at \(p^\ell\) requires
at least one bubble center converging to \(p^\ell\). Since \(p^1,\ldots,p^I\)
are distinct and each nonflat bubble satisfies \(\operatorname{ind}(V^i)\ge1\),
the index bound \(\sum_{i=1}^N\operatorname{ind}(V^i)\le I\) forces \(N=I\) and
\(\operatorname{ind}(V^i)=1\) for each \(i\).$\square$
\end{proof}

\vskip0.33in

\section{\bf Model bubble energies and seven-dimensional counterexamples}\label{sec:examples}

The radial exterior-ball solution has Morse index one by separation of
variables \cite[Remark~1.3]{BKPeriodic}. We include the calculation, compute
its exact Hessian energy, and illustrate two rescaling regimes. We also
show that the uniform bounded-index Hessian energy bound fails in
dimension seven.

\vskip0.23in

\begin{proposition}\label{prop:exterior-ball-energy}
For every \(n\ge3\), the unit exterior-ball solution
\(V(x)=(n-2)^{-1}(1-|x|^{2-n})\) on \(\{|x|>1\}\), extended by zero to
\(\overline B_1\), satisfies
\[
  \operatorname{ind}(V)=1,\qquad
  E_n^{\mathrm{ball}}:=\mu_V(\R^n)
  =|\mathbb S^{n-1}|\,[n(n-1)]^{n/2-1},
\]
and both quantities are invariant under rigid motions and Bernoulli rescaling.
\end{proposition}

\begin{proof}
At \(r=|x|>1\), \(D^2V(x)\) has radial eigenvalue \(-(n-1)r^{-n}\) and \(n-1\)
tangential eigenvalues \(r^{-n}\), so \(|D^2V(x)|^2=n(n-1)r^{-2n}\). Integrating
in polar coordinates gives
\[
  \mu_V(\R^n)
  =|\mathbb S^{n-1}|\,[n(n-1)]^{n/2}\int_1^\infty r^{n-1-n^2}\dd r
  =|\mathbb S^{n-1}|\,[n(n-1)]^{n/2-1},
\]
since \(\int_1^\infty r^{n-1-n^2}\dd r=(n^2-n)^{-1}\).

\vskip0.23in

To compute \(\operatorname{ind}(V)\), recall \(H=n-1\) on \(\partial B_1\) and
expand the boundary trace \(\phi|_{\partial B_1}\) of an admissible test
function \(\phi\in C_c^\infty(\R^n)\) in an \(L^2(\mathbb S^{n-1})\)-orthonormal
basis of spherical harmonics \(\{Y_{\ell m}\}\) of degree \(\ell\ge0\) with
coefficients \(c_{\ell m}\). By Dirichlet's principle on \(\R^n\setminus\overline{B_1}\),
\(\int_{|x|>1}|\nabla\phi|^2\dd x\) is bounded below by the Dirichlet energy of
the decaying harmonic extension \(\sum_{\ell,m}c_{\ell m}r^{-(\ell+n-2)}Y_{\ell m}\),
giving
\[
  \mathcal Q_V[\phi;\R^n]
  \ge \sum_{\ell,m}(\ell+n-2)|c_{\ell m}|^2
      -(n-1)\sum_{\ell,m}|c_{\ell m}|^2
  =\sum_{\ell,m}(\ell-1)|c_{\ell m}|^2.
\]
Hence \(\mathcal Q_V\) is nonnegative on the codimension-one subspace of traces
with zero spherical average (\(c_{00}=0\)), giving \(\operatorname{ind}(V)\le1\).
Conversely, take \(\chi\in C_c^\infty([0,\infty))\), equal to one on
\([0,1]\) and zero on \([2,\infty)\). The test
\(\phi_R(x)=\chi(|x|/R)|x|^{2-n}\) on \(\{|x|\ge1\}\) has a smooth
compactly supported extension to \(\R^n\), obtained by changing the formula
only near the origin. The uncut radial mode has form
\(-|\mathbb S^{n-1}|\), and the tail and cutoff change its Dirichlet energy
by \(O(R^{2-n})\). Thus \(\mathcal Q_V[\phi_R;\R^n]<0\) for large \(R\),
giving \(\operatorname{ind}(V)=1\). Rigid motions preserve both quantities.
For the scaled profile \(\lambda V((x-z)/\lambda)\), the pullback of
tests multiplies the quadratic form by \(\lambda^{n-2}\), and the
Hessian energy is invariant.$\square$
\end{proof}

\vskip0.23in

\begin{remark}\label{rem:contrasting-examples}
In Theorem~\ref{thm:local-bubble-tree}, if \(m_p\) exterior-ball bubbles
concentrate at \(p\), they contribute atomic mass \(m_p E_n^{\mathrm{ball}}\)
and consume \(m_p\) units of Morse index. Even for a single shrinking exterior
ball \(u_\lambda(x)=\lambda V(x/\lambda)\) (\(\lambda\downarrow0\)), one sees
both hallmarks of Proposition~\ref{prop:intermediate-host}: \(u_\lambda\to 0\)
uniformly, \(\one_{\{u_\lambda>0\}}\to 1\) in \(L^1_{\mathrm{loc}}(\R^n)\), and
\(\mu_{u_\lambda}\rightharpoonup E_n^{\mathrm{ball}}\delta_0\) with
\(\operatorname{ind}(u_\lambda;B_2)=1\) for small \(\lambda\). By contrast, the
collapsing stable two-plane family
\(u_\lambda(x)=(x_n-\lambda)_++(-x_n-\lambda)_+\) converges uniformly to the
centered unit wedge \(|x_n|\) with \(\one_{\{u_\lambda>0\}}\to 1\) in
\(L^1_{\mathrm{loc}}(\R^n)\), while carrying identically zero Hessian energy
\(\mu_{u_\lambda}\equiv0\) at every \(\lambda>0\).
\end{remark}

\vskip0.23in

In dimension seven, a smooth minimizing leaf with nonflat conical blow-down
gives a family of Morse index zero whose interior Hessian energies grow at
least logarithmically. Hence the uniform local mass bound underlying
Theorem~\ref{thm:local-bubble-tree} fails in dimension seven.

\vskip0.23in

\begin{proposition}\label{prop:seven-energy}
In \(\R^7\), there exists a one-parameter family of entire classical minimizers
\(\{u_\lambda\}_{0<\lambda<1/2}\) satisfying \(0\in\FB(u_\lambda)\) and
\(\operatorname{ind}(u_\lambda;B_2)=0\) such that
\[
  \mu_{u_\lambda}(B_1)\ge c\log(1/\lambda)-C
  \qquad\forall\,\lambda\in(0,1/2),
\]
for positive constants \(c,C>0\) independent of \(\lambda\).
\end{proposition}

\begin{proof}
Let \(C_0\) be the seven-dimensional nonflat minimizing cone of
De Silva and Jerison \cite{DSJ}, and let \(V\) be a smooth entire minimizing leaf
with blow-down \(C_0\) provided by \cite[Theorem~1.5]{ESV}. Because \(C_0\) is
\(1\)-homogeneous and nonflat, its positive-phase Hessian does not vanish
identically on \(A(1,2)\cap\{C_0>0\}\), so there exists a compact set
\(K\Subset A(1,2)\cap\{C_0>0\}\) on which \(\int_K|D^2C_0|^7\dd x>0\). Since
\(R^{-1}V(Rx)\to C_0(x)\) locally uniformly as \(R\to\infty\), interior
harmonic estimates yield smooth convergence on \(K\), whence there exist
\(R_0\ge1\) and \(c_0>0\) such that
\[
  \mu_V(B_{2R}\setminus B_R)
  \ge \int_K\bigl|D^2\bigl(R^{-1}V(Rx)\bigr)\bigr|^7\dd x
  \ge c_0>0
  \qquad\forall\,R\ge R_0.
\]
Fix any point \(z_0\in\FB(V)\) and define the rescaled entire minimizers
\(u_\lambda(x)=\lambda V(z_0+x/\lambda)\). By scale invariance,
\(\mu_{u_\lambda}(B_1)=\mu_V(B_{1/\lambda}(z_0))\). For small \(\lambda>0\),
the ball \(B_{1/\lambda}(z_0)\) contains \(\asymp\log(1/\lambda)\) pairwise
disjoint dyadic annuli \(B_{2^{k+1}R_0}(0)\setminus B_{2^k R_0}(0)\), each
contributing at least \(c_0\) to \(\mu_V\). This proves the lower bound
for small \(\lambda\), and increasing \(C\) covers the rest of \((0,1/2)\).
The functions \(u_\lambda\) are minimizers by scaling and translation,
so their Morse index on \(B_2\) is zero.$\square$
\end{proof}

\vskip0.23in

\begin{remark}\label{rem:seven-outlook}
In dimension seven, Proposition~\ref{prop:seven-energy} exhibits conical
annuli that retain positive Hessian energy on each dyadic scale even at
Morse index zero.
\end{remark}

\vskip0.63in

\section*{Declarations}
\noindent\textbf{Availability of data and materials.} No datasets were
generated or analysed during the current study.

\medskip
\noindent\textbf{Conflict of interest.} The authors declare no conflict
of interest.

\medskip
\noindent\textbf{Use of artificial intelligence.} AI tools were used to assist in the exploration and verification of mathematical proofs, as well as the translation and linguistic polishing of the manuscript. The authors
take full responsibility for the mathematical content and the final manuscript.

\vskip0.63in

\end{document}